\documentclass[11pt, reqno]{amsart}

\usepackage[utf8]{inputenc}
\usepackage[T1]{fontenc}
\usepackage[english]{babel}
\usepackage[margin=2.5cm, includefoot, footskip=30pt, voffset=.1in, headsep=.25in, textheight=625pt]{geometry}
\usepackage{amsfonts, amssymb, amsmath} 
\usepackage{color}
\usepackage{graphicx}
\usepackage[all, cmtip]{xy}
\usepackage{eurosym}
\usepackage{amsthm}
\usepackage{stmaryrd}
\usepackage{appendix}
\usepackage{scrextend}
\usepackage{lipsum}
\usepackage[linkcolor=blue, citecolor=green, colorlinks=true]{hyperref}
\usepackage{array}
\usepackage{indentfirst}
\usepackage[Lenny]{fncychap}

\usepackage{CJKutf8}
\usepackage{ucs}

\usepackage{fancyhdr}
\makeatletter
\let\@oddfoot\@empty
\let\@evenfoot\@empty
\makeatother

\usepackage{libertine} \usepackage[libertine]{newtxmath}
\usepackage{chngcntr}
\counterwithout{equation}{section}
\numberwithin{equation}{section}	

\newcommand{\R}{\mathbf{R}}
\newcommand{\C}{\mathbf{C}}

\newcommand{\N}{\mathbf{N}}
\newcommand{\A}{\mathbf{A}}
\newcommand{\h}{\mathcal{H}}

\renewcommand{\P}{\mathbf{P}}

\newcommand{\Q}{\mathbf{Q}}
\newcommand{\SL}{\mathrm{SL}}
\newcommand{\GL}{\mathrm{GL}}

\newcommand{\PGL}{\mathrm{PGL}}

\renewcommand{\Re}{\mathrm{Re}}

\newcommand{\vol}{\mathrm{vol}}
\newcommand{\id}{\mathrm{id}}

\newcommand{\tr}{\mathrm{tr} \ }
\newcommand{\itex}{\item[$\bullet$]}
\newcommand{\tq}{\:: \:}

\newcommand{\p}{\mathfrak{p}}
\newcommand{\q}{\mathfrak{p}}
\newcommand{\dd}{\mathrm{d}}

\renewcommand{\det}{\mathrm{det}}

\renewcommand{\div}{\: | \:}
\renewcommand{\mod}{\mathrm{ mod } \ }

\newcommand{\spec}{\mathrm{spec}}

\newcommand{\geom}{\mathrm{geom}}
\renewcommand{\tr}{\mathrm{tr }}
\renewcommand{\q}{\mathfrak{q}}
\renewcommand{\d}{\mathfrak{d}}
\renewcommand{\c}{\mathfrak{c}}
\renewcommand{\l}{\mathfrak{l}}
\renewcommand{\p}{\mathfrak{p}}

\newcommand{\f}{\mathfrak{f}}
\newcommand{\n}{\mathfrak{n}}
\newcommand{\m}{\mathfrak{m}}
\renewcommand{\dim}{\mathrm{dim } \ }
\renewcommand{\vol}{\mathrm{vol }}

\renewcommand{\GL}{\mathrm{GL}}		
\renewcommand{\SL}{\mathrm{SL}}

\renewcommand{\PGL}{\mathrm{PGL}}
\renewcommand{\O}{\mathcal{O}}

\renewcommand{\ell}{\mathrm{ell}}
\newcommand{\Pl}{\mathrm{Pl}}
\newcommand{\JL}{\mathrm{JL}}
\newcommand{\PPP}{\phi}
\renewcommand{\h}{\mathfrak{h}}

\newcommand{\rk}{\noindent \textit{Remark. }}
\newcommand{\rks}{\noindent \textit{Remarks. }}

\theoremstyle{plain}
\newtheorem{theorem}{Theorem}
\newtheorem{lemma}{Lemma}
\newtheorem{proposition}{Proposition}[section]

\newtheorem{corollary}[proposition]{Corollary}
\theoremstyle{definition}

\newcolumntype{L}{>{$}c<{$}}

\begin{document}

\title{Counting and Equidistribution for Quaternion Algebras}

\author{Didier Lesesvre}

\address{
              Department of Mathematics ,
              Sun Yat-Sen University,
              Zhuhai, Guangdong, China
              }

\begin{abstract}
We aim at studying automorphic forms of bounded analytic conductor in the division quaternion algebra setting. We prove the equidistribution of the universal family with respect to an explicit and geometrically meaningful measure. It leads to answering the Sato-Tate conjectures in this case, and contains the counting law of the universal family, with a power savings error term in the totally definite case.
\end{abstract}

\maketitle

\section{Introduction}
\label{sec:intro}

\subsection{Landscape}
\label{subsec:landscape}

Let $F$ be a number field of degree $d$ over $\Q$. Let $\A$ denote the ring of adeles of $F$. We consider a division quaternion algebra $B$ over $F$, and write $R$ for the places of $F$ where $B$ is not split. We introduce the group of projective units $G = Z \backslash B^\times$, where $Z$ denotes the center of $B^\times$. Let $\mathcal{A}(G)$ denote the universal family of $G$, that is the set of all irreducible automorphic infinite dimensional representations of $G(\A)$. A deep understanding of $\mathcal{A}(G)$ is of fundamental importance in the theory of automorphic forms. 

In order to determine its actual size and some sharper statistical properties, as densities or equidistribution, we need to truncate it for we then deal with a finite set, hence we need a suitable notion of size to do so. Turn for a moment to a more usual setting: the one of general linear groups. The universal family $\mathcal{A}(G)$ embeds, via the Jacquet-Langlands correspondence, as a subfamily of the universal family $\mathcal{A}(\PGL(2))$, composed of all the cuspidal automorphic representations of $\PGL(2)$. In the latter context, even in the broader setting of cusp forms on general linear groups, \cite{iwaniec_perspectives_2000} have defined a good notion of size, given by the analytic conductor. It is a positive real number $c(\pi)$\label{c} defined from the functional equation satisfied by the finite part L-function $L(s,\pi)$\label{L} associated to $\pi\in \mathcal{A}(\PGL(2))$, which takes the form
\begin{equation}
L(1-s, {\pi}) = \bar{\varepsilon}_\pi X(s, \pi) L(s, \pi), 
\end{equation}

\noindent where $\varepsilon_\pi$ is the root number\label{epsilon} of $\pi$. The completing factor $X(s, \pi)$\label{XL} takes value $\varepsilon_\pi$ at the central point $\frac{1}{2}$, and the analytic conductor is defined to be $c(\pi) = |X'(\frac{1}{2}, \pi)|$ following the presentation of  \cite{conrey_integral_2005}. The function $X(s, \pi)$ involves the usual arithmetic conductor as well as archimedean gamma factors, so that the analytic conductor encapsulates the complexity of $\pi$. It allows to truncate the universal family of $\PGL(2)$, and hence the one of $G$, to a finite set \cite{brumley_effective_2006}. The truncated universal family may then be introduced as
\begin{equation}
\label{UF}
\mathcal{A}(Q) = \{\pi \in \mathcal{A}(G) \ : \   c(\pi) \leqslant Q\}, \qquad Q \geqslant 1.
\end{equation}

The problem of counting automorphic representations ordered by analytic conductor goes back to the work of  \cite{iwaniec_perspectives_2000}.  In this article, we seek to prove certain basic properties of this family, such as determining its asymptotic growth, establishing global equidistribution with respect to a geometric significant measure, and proving the validity of the Sato-Tate conjecture in this setting.

\subsection{Analogy with the height on an algebraic variety}
\label{subsec:analogy}

The counting problem admits an interesting analogy with the well-known question of counting rational points of bounded height on a smooth projective variety over a number field. The absolute Weil height is the proper notion of size in this setting and is defined by
\begin{equation}
\label{Weil-height}
h_{\mathbf{P}^n}(x) = \prod_v \max_{0 \leqslant i \leqslant n} |x_i|_v^{1/[F:\Q]}, \qquad x=(x_i)_{0 \leqslant i \leqslant n} \in \P^n(F).
\end{equation}

\noindent where the product runs over the places of $F$ and does not depend on the choice of homogeneous coordinates. Given any projective variety $V$ over $F$ endowed with a fixed embedding $\iota$ into the projective space $\P^n(F)$, a height function on $V$ can be defined by pulling back the Weil height on $\P^n(F)$, setting
\begin{equation}
\label{pullback}
h_V(x) = h_{\P^n}(\iota(x)), \qquad x \in V.
\end{equation}

The most natural setting for considering such generalized questions is the one of Fano varieties, where there are precise conjectures due to \cite{batyrev_sur_1990, peyre_hauteurs_1995}. On those grounds,  \cite{northcott_periodic_1950} proved the finiteness of the set of points of bounded height for the projective spaces, refined by \cite{schanuel_heights_1964} in an asymptotic counting law.

\begin{theorem}[Schanuel]
\index{Schanuel theorem}
\label{thm:schanuel}
For all $n \geqslant 1$, there exists $C_n > 0$ such that for any $Q \geqslant 1$, 
\begin{equation*}
\# \left\{ x \in \P^n(F) \ : \ h(x) \leqslant Q \right\} = C_nQ^{n+1} + 
\left\{
\begin{array}{cl}
O\left(Q \log Q\right) & \text{ if } n=1, \ F = \Q ; \\
O\left(Q^{n-1/[F:\Q]} \right) & \text{otherwise.}
\end{array}
\right. 
\end{equation*}
\end{theorem}

The analogy between the Schanuel theorem on counting rational points on projective spaces and the problem of counting automorphic cusp forms on $\GL(n)$ has been particularly stressed recently, according to Sarnak. The case of quaternion algebras can be embedded in $\GL(2)$ so that, following the above analogy, the notion of analytic conductor we use in our main theorem is inspired by the procedure for heights: given the by now standard notion of analytic conductor for $\GL(2)$, we pull it back to automorphic forms on quaternion algebras via the associated identity map between their dual groups, hence defining the notion of analytic conductor in our setting.

\subsection{Counting law for the universal family}
\label{subsec:counting}

The first result of this article gives an asymptotic formula for the cardinality
\begin{equation}
N(Q) = \# \mathcal{A}(Q), \qquad Q \geqslant 1,
\end{equation}

\noindent  Petrow recently handled the problem in a fairly general fashion for automorphic forms on tori \cite{petrow_weyl_2018}. The case of the universal family for $\GL(2)$ is handled in the preprint \cite{brumley_counting_2018}. For division quaternion algebras, the counting law is provided by the following statement.

\begin{theorem}[Counting law for quaternion algebras]\label{thm-count}
There exists $C > 0$ such that for any $Q \geqslant 1$, 
\begin{equation*}
N(Q) =  C Q^2 + 
\left\{
\begin{array}{cl}
O\left(Q^{1+\varepsilon}\right) & \text{ if } F = \Q \text{ and $B$ totally definite, for all } \varepsilon >0; \\
O\left(Q^{2-\delta_F} \right) & \text{ if $F \neq \Q$} \text{ and $B$ totally definite} ;\\
O\left(\frac{Q^2}{\log Q} \right) & \text{ if $B$ is not totally definite}.
\end{array}
\right.
\end{equation*}

\noindent The constant $C>0$ is defined explicitly in \eqref{constant}, and $\delta_F = 2({1+[F:\Q]})^{-1}$.
\label{thm1}
\label{deltaF}
\end{theorem}

\noindent \textit{Remarks.} The form of this asymptotic growth appeals some comments.
\begin{itemize}
\item[(i)] There is a similarity between the error term in Theorem \ref{thm-count} and that of the classical result of Schanuel in Theorem \ref{thm:schanuel} on the number of rational points of bounded height in projective spaces. His result, when specialized to $F=\Q$, has an error term that picks up an additional small quantity, namely a power of log, to be compared to the $Q^\varepsilon$ of Theorem \ref{thm-count}.
\item[(ii)] The presence of a power savings error term in the totally definite case, \textit{i.e.} when every archimedean place is ramified, is noteworthy. This feature is lost without this assumption, like the corresponding result \cite{brumley_counting_2018} for $\GL(2)$, where only a logarithmic savings is obtained. The reason for this difference lies in the passage from smooth to sharp counting at archimedean places, see Section \ref{subsec:smoothing-error}.
\item[(iii)] The assumption that $B$ is a division quaternion algebra induces an automorphic compact quotient, hence avoids considerations due to the continuous part of the automorphic spectrum, see Section \ref{sec:trace-formula}.
\item[(iv)] The center has been removed for technical purposes and to avoid to deal with the central terms in the Selberg trace formula. All the methods are expected to carry on to a setting considering the center without considerable adaptation.
\end{itemize}

\medskip

The precise knowledge of the constant $C$ unveils a lot of information, and its geometrical interpretation has considerable importance as in the conjectures of Peyre. An explicit and meaningful formulation of the constant is given below, in the context of the equidistribution properties of $\mathcal{A}(G)$, and shows striking similarities with the ones computed for algebraic varieties \cite{chambert-loir_igusa_2010}.

\subsection{Equidistribution of the universal family}
\label{subsec:equid}

Beyond estimating the size of the universal family lies the question of the geometrical distribution of the automorphic representations of $G$. A good formulation of the problem is developed in the work of \cite{sarnak_families_2016} and is to find a measure with respect to which the universal family equidistributes\index{equidistribution}, what is carried on in this section after giving a glance at the topological and measurable structure the universal family is endowed with. 

Each local unitary dual group $\widehat{G}_v$ is endowed with the Fell topology and the product $\prod_v \widehat{G}_v$ is then given the product topology. Introduce the measure $\mu$ on $\prod_v \widehat{G}_v$ that assigns to every basic open set $X=\prod_v X_v$, \textit{i.e.} where $X_v$ is an open set of $\widehat{G}_v$ and $X_v=\widehat{G}_v$ for all but finitely many $v$, the positive real number
\begin{equation}
\mu(X)= \int_{X}^\star \frac{\dd\pi}{c(\pi)^2},
\label{mu}
\end{equation}

\noindent where the regularized integral\index{regularized integral} is defined as
\begin{equation}
\label{regularized-integral-mu}
\zeta^{\star}(1) \prod_v \zeta_v(1)^{-1}  \int_{{X_v}} \frac{\dd\pi_v}{c(\pi_v)^2}.
\end{equation}

Here $\zeta_v$\label{zeta} is the local zeta function associated to $F_v$, the notation $\zeta^{\star}(1)$ stands for the residue of the Dedekind zeta function of $F$ at $1$, and $\dd\pi_v$ is the Plancherel measure on $\widehat{G}_v$, introduced and normalized according to the convention in Section \ref{sec:plancherel-formulas}. 

\medskip

\rks This integral is not as disturbing as it seems for the following reasons.
\begin{itemize}
\item[(i)] The Plancherel measure is supported on the tempered dual; since tempered representations are generic, the conductors appearing in the integral are well-defined for the sets actually arising in what follows, see Section \ref{sec:conductors}.
\item[(ii)] It is by no mean obvious that the integral \eqref{mu} actually converges. It is the case and Section \ref{subsec:constant} contains the explicit computations of the local factors ensuring the convergence as well as motivating the regularization, see Section \ref{sec:integralscv}.
\end{itemize}

\medskip

The measure $\mu$ has finite total mass $\|\mu\|$. All the definitions are now in place to uncover the expression of the leading constant in Theorem \ref{thm-count}, namely
\begin{equation}
C =  \frac{1}{2} \mathrm{vol}(G(F) \backslash G(\A)) \|\mu\|,
\label{constant}
\end{equation}
\noindent where the measure giving the volume of the automorphic quotient $G(F) \backslash G(\A)$ is normalized as in Section \ref{subsec:Haar-measures}. The main result is the following one.

\begin{theorem}[Equidistribution for quaternion algebras]
\label{thm-equid}
The universal family of $G$ equidistributes with respect to the measure $\mu$, in the following sense. For every relatively quasi-compact open set $X$ of $\prod_v \widehat{G}_v$ with boundary of measure zero,
\begin{equation}
\frac{\#\{\pi \in \mathcal{A}(Q): \pi\in X\}}{N(Q)} \longrightarrow \frac{\mu}{\|\mu\|}(X), \qquad \text{as} \quad Q\rightarrow \infty.
\end{equation}
\end{theorem}


Once this global equidistribution result stated, the Sato-Tate conjecture questions the behavior of the projections $\mu_\p$ of the limit measure on the local components $\widehat{G}_\p$ when the norm of $\p$ grows. Let $T_c$ be the subgroup of diagonal matrices in $\mathrm{SU}(2)$ and $W$ the associated Weyl group. On the common ground where all the representations in the support of the Plancherel measures of $G_\p$ live, given by the tempered Satake parameters space $T_c/W$, the Sato-Tate question acquires a precise meaning and local representations are equidistributed with respect to the half-circle measure.

\begin{corollary}[Sato-Tate for quaternion algebras]
\label{coro:ST}
For all $\phi \in C(T_c/W)$, 
\begin{equation}
\int_{T_c/W} \widehat{\phi}(x) \dd\mu_\p(x) \longrightarrow \int_{T_c/W} \widehat{\phi}(x) \dd\mu^{\mathrm{ST}}(x), \qquad \text{as} \quad N\p \longrightarrow \infty, 
\end{equation}
\noindent where $\dd\mu^{\mathrm{ST}}$ is the Sato-Tate measure on the half-circle, \textit{i.e.}
\begin{equation}
\dd\mu^{\mathrm{ST}}(x) = \frac{1}{\pi} \sqrt{1-\frac{x^2}{4}}\dd x.
\end{equation}
\end{corollary}

\subsection{Organization of this article}
\label{subsec:organisation}

Section \ref{sec:prelim} is mainly devoted to introducing notations, stating the precise definition of the analytic conductor, and fixing the normalizations of measures. We recall some facts about equidistribution and spectral tools required to reduce Theorem \ref{thm-equid} to a statement amenable to trace formula methods, among which are the Sauvageot density theorem and Paley-Wiener theorems. Section \ref{sec:combi} provides a decomposition of the universal family into harmonic subfamilies obtained by fixing certain spectral data. In Section \ref{sec:fts} we show that the proportion of automorphic forms we seek to estimate can be expressed as a spectral side of the Selberg trace formula, hence can be expressed in terms of orbital integrals. The main asymptotic term involved in this proportion comes from the contribution of the identity, which we evaluate in Section \ref{sec:identity}. Other spectral and geometric terms arise in the trace formula. The spectral ones are those coming from undesired characters, side effects of the smoothing procedure for test functions at archimedean places and  the contribution of the complementary spectrum: they are precisely bounded in Section \ref{sec:REMAINDER}. The geometric ones are those coming from the orbital integrals associated to other terms than the identity, and they are bounded in Section \ref{sec:elliptic}, opening the path to the claimed asymptotic development. The ultimate Section \ref{sec:ST} builds on the known Plancherel measures in the split case in order to prove that the limit measure with respect to which the universal family equidistributes satisfies the Sato-Tate equidistribution conjecture.

\subsection{acknowledgements}
I am infinitely indebted to my advisor, Farrell Brumley, for having trusted in me for handling this problem and for having been strongly present and implied in its resolution. I am grateful to the time Gergely Harcos and Philippe Michel granted me in carefully reading the thesis from which this article blossomed. I would like to thank Valentin Blomer, Andrew Corbett, Mikolaj Fraczyk, Élie Goudout and Guy Henniart for many enlightening discussions. This work has been faithfully supported by the ANR 14-CE25 PerCoLaTor, the Fondation Sciences Mathématiques de Paris and the Deutscher Akademischer Austauschdienst.  At last, nothing would have been done without the warm and peaceful environment of the institutions which hosted me during these years: Université Paris 13, 	École Polytechnique Fédérale de Lausanne, Georg-August Universität and Sun Yat-Sen University.

\section{Groundwork}
\label{sec:prelim}

We denote by $v$ the places of $F$, $\p$ the non-archimedian ones, and $\O_\p$ the ring of integers of $F_\p$ for a finite place $\p$. The finite set $R$ of ramification places of $B$ determines it up to isomorphism. From now on, Latin letters $q, d, m$, etc. will denote usual integers, while Gothic letters $\q, \d, \mathfrak{m}$, etc. will denote ideals of integer rings. 

\subsection{Analytic conductor}
\label{subsec:conductor}
\label{sec:conductors}

In order to make sense of the problem, we need to define precisely the notion of size we choose for representations. It is the analytic conductor, which we introduce in this section. We will work with $B^\times$ more than with $G$, for it lightens notations. This local convention makes no harm, for we view a representation $\pi$ of $G(\A) = PB^{\times}(\A)$ as a representation of $B^\times(\A)$ with trivial central character. By Flath's theorem, an irreducible admissible representation of $B^\times(\A)$ decomposes in a unique way as a restricted tensor product $\pi = \otimes_v \pi_v$ of irreducible smooth representations where almost every component $\pi_v$ is unramified. We want first to define the conductor for the local components $\pi_v$. 

\medskip

The Jacquet-Langlands correspondence allows to reduce to the $\GL(2)$ case, and in this one only infinite-dimensional representation arise. Indeed, since the universal family excludes global characters, a representation $\pi$ in it is generic. The Jacquet-Langlands correspondence preserves genericity hence, as shown on the diagram below, associates to $\pi$ a generic representation $\JL(\pi)$ of $\GL(2)$, thus also its local components $\JL(\pi)_v$. These local components are also the images by the local Jacquet-Langlands correspondence $\JL_v(\pi_v)$ of the local components of $\pi$. 
$$
\xymatrix{
\, \substack{\pi \in \mathcal{A}({B^\times})} \quad \ar[rr]^{\mathrm{JL}} \ar[d]_v & & \quad \, \substack{\mathrm{JL}(\pi) \in \mathcal{A}(\GL(2)) \\ \mathrm{generic}} \ar[d]^v \\
\pi_v \quad \ar[rr]^{\mathrm{JL}_v}_{\id \ \mathrm{if} \ v \notin R} & & \quad \substack{\mathrm{JL}(\pi)_v \\ \mathrm{generic}}
}
$$

At split places, the local Jacquet-Langlands correspondence is the identity, for then $B_\p^\times \simeq \GL(2, F_\p)$. The correspondence is unique, thus the local components at split places $\pi_v$ are generic hence infinite-dimensional, proving the claim. 

\subsubsection{Non-archimedian case}

For finite split places $\p$, by definition $B_{\p} \simeq M\left(2, F_{\p}\right)$ so that $B^\times_{\p}$ is isomorphic to $\GL(2, F_\p)$. The notion of local conductor for irreducible smooth infinite-dimensional representations of $\GL(2)$ has been introduced by \cite{casselman_results_1973}. Consider the sequence of compact open congruence subgroups
\begin{equation}
\label{filtration}
K_{0, \p}\left(\p^r\right) = 
\left\{
g \in \GL\left(2, \O_{\p}\right) \tq 
g \equiv
\left(
\begin{array}{cc}
\star & \star \\
0 & \star
\end{array}
\right)
\mod \p^r
\right\} \subseteq B^\times_{\p}, \qquad r \geqslant 0.
\end{equation}

The multiplicative and analytic conductors of an  irreducible admissible infinite-dimensional representation $\pi_\p$ of $B^\times_\p$ with trivial central character are then respectively defined by
\begin{equation}
\label{conductor-XXX}
\c(\pi_\p) = \p^{\f(\pi_{\p})} \quad \text{and} \quad  c(\pi_\p) = N\c(\pi_\p),
\end{equation}

\noindent where 
\begin{equation}
 \f(\pi_{\p}) = \min \left\{r \in \N \tq \pi_{\p}^{K_{0, \p}\left(\p^r\right)} \neq 0\right\}.
\end{equation}

\noindent The existence of the conductor is guaranteed by the work of  \cite{casselman_results_1973}, who also states that the growth of the dimensions of the fixed vector spaces are given by
\begin{equation}
\dim \pi_\p^{K_{0, \p}\left(\p^{\f(\pi_\p)+i}\right)} = i+1, \qquad i \geqslant 0.
\label{one}
\end{equation}

\subsubsection{Archimedian case}

The archimedian part of the conductor is introduced by  \cite{iwaniec_perspectives_2000}. It is built on the archimedean factors completing the L-functions associated to automorphic representations. The archimedean L-factors are of the form
\begin{equation}
L(s, \pi_v) = \prod_{j=1}^2 \Gamma_v(s-\mu_{j,\pi}(v)),
\end{equation}

\noindent where $\Gamma_v(s) = \pi^{-s/2}\Gamma(s/2)$ and the $\mu_{j,\pi}(v)$ are complex numbers. The analytic conductor is then locally defined to be, for $v \div \infty$, 
\begin{equation}
c_v(\pi) = \prod_{j=1}^2 \left(1 + |\mu_{j, \pi}(v)|\right).
\end{equation}

\begin{rk}
We cannot avoid, following Iwaniec and Sarnak, this archimedean part of the conductor and only consider its arithmetic component. Indeed, we aim at counting irreducible admissible infinite-dimensional representations of bounded conductor, but this family would be infinite without control of the archimedian conductor. For instance the family of modular forms of level one and arbitrary weights constitutes an infinite family of fixed arithmetic conductor: they give rise to the discrete series.
\end{rk}

\subsubsection{Non-split case}

Via the Jacquet-Langlands correspondence, the non-split case is reduced to the already treated split one, analogously with the pullback of heights for algebraic varieties. The conductor of an irreducible admissible representation $\pi_{v}$ of $B^\times_{v}$ is defined as the conductor of its Jacquet-Langlands transfer
\begin{equation}
c\left(\pi_{v}\right) = c\left(\JL\left(\pi_{v}\right)\right).
\end{equation}

\subsubsection{Characters}

For now conductors have been defined only for generic representations. However, characters can arise as local components at ramified places as discussed above. Every character of $B_\p^\times$ is a composition
\begin{equation}
B_\p^\times \longrightarrow F_\p^\times \longrightarrow \C,
\end{equation}

\noindent where the first application is the reduced norm, and the second one a character of $F_\p^\times$. In other words, every character of $B_\p^\times$ is of the form $\chi_0 \circ N$ where $\chi_0$ is a character of $F_\p^\times$ and $N$ the reduced norm on $B_\p^\times$. In order to stay consistent, define the conductor of a local character at a ramified place as the conductor of its Jacquet-Langlands embedding in $\GL(2)$, which is defined based on the associated functional equation. Since the character $\chi_0 \circ N$ is sent to the twisted Steinberg representation $\mathrm{St} \otimes \chi_0$, it follows explicitly
\begin{equation}
\c(\chi_0 \circ N) = \left\{
\begin{array}{cl}
\p & \text{if $\chi_0$ unramified;} \\
\c(\chi_0)^2 & \text{if $\chi_0$ ramified.}
\end{array}
\right.
\end{equation}

\subsubsection{Global analytic conductor}

We introduce for an irreducible admissible representation of $B^\times(\A)$ decomposed into $\pi = \otimes_{v} \pi_{v}$ its global analytic conductor
\begin{equation}
c\left(\pi\right) = \prod_{v} c\left(\pi_v\right).
\end{equation}

\noindent This gives a well-defined notion of conductor, for the $\pi_v$ are almost everywhere unramified, thus of conductor one. It extends to a definition for representations of $G(\A)$, viewed as automorphic representations of $B^\times(\A)$ with trivial central characters.

\medskip

\noindent \textit{Remarks.}  Analogously to what happens for general linear groups, the conductor could have been defined directly from the $L$-functions associated to automorphic representations of quaternion algebras. These are provided by the \cite{godement_zeta_1972} construction and would avoid the appeal to an embedding in $\mathrm{GL}(n)$. The Jacquet-Langlands correspondence preserves the notion of $L$-function and hence also makes this notion of conductor for $G$ compatible with the one obtained by pulling back the conductor on $\GL(2)$, thus this choice of exposition makes no harm compared to directly defining the conductor from the associated $L$-functions on $G$. Thus both choices of definition of the conductor coincide.

\subsection{Normalization of measures}
\label{subsec:measures}
\label{sec:plancherel-formulas}
\label{subsec:Haar-measures}

At the non-archimedean places, the measure taken on $G_\p$ is the Haar measure $\mu_\p$ normalized so that $K_\p = \PGL(2,\O_{\p})$, in the split case, or $K_\p = \mathfrak{o}_{\p}^\times$ the units of a maximal order of $B_\p$, in the non-split case, gets measure one. This normalisation is independent of the chosen maximal order \cite{hull}. For the archimedean places, we choose the Haar measure normalized so that the maximal compact subgroup gets measure one. 

We now turn to the associated local dual groups. Denote $\mathcal{H}(G_v)$ the Hecke algebra of $G_v$, that is the algebra consisting of compactly supported complex-valued functions on $G_v$, locally constant at finite places, smooth at archimedian ones. Let $\mathcal{H}(G(\A))$ be the Hecke algebra of $G(\A)$. It is the algebra generated by the restricted products $\phi = \prod_v \phi_v$, where $\phi_v$ is a function of $\mathcal{H}(G_v)$ and almost every local component $\phi_\p$ is equal to $\mathbf{1}_{K_{\p}}$. For such a function $\phi \in \mathcal{H}(G(\A))$, we extend the action of $\pi$ to $\mathcal{H}(G(\A))$, $\pi(\phi)$ acting by the mean action of $\pi$ over $G$ with weight $\phi$, that is to say
\begin{equation}
\pi(\phi) = \int_{G(\A)} \phi(g) \pi(g) \dd g.
\end{equation}

\noindent This defines a Hilbert-Schmidt integral operator of trace class, thus we can define its Fourier transform by
\begin{equation}
\widehat{\phi}(\pi) = \tr \ \pi(\phi) = \tr \left( v \mapsto \int_{G} \phi(g)\pi(g)v \dd g \right).
\end{equation} 

The unitary dual group $\widehat{G}_v$ is endowed with its usual Fell topology and Plancherel measure associated with the measure chosen on $G_v$: it is the unique positive Radon measure $\mu^{\Pl}_v$ on $\widehat{G}_v$ such that the Plancherel inversion formula of Harish-Chandra holds, \textit{i.e.} for functions $\phi_v$ in the Hecke algebra $\mathcal{H}(G_v)$, we have
\begin{equation}
\int_{\widehat{G}_v} \widehat{\phi}_v(\pi_v) \dd\mu_v^{\Pl}(\pi_v) = \phi_v(1).
\label{plancherel}
\end{equation}

\noindent From now on, integrals on $\widehat{G}_v$ will be written with the convention that $\dd\pi_v = \dd\mu_v^{\Pl}(\pi_v)$, leading to no ambiguity. On $\widehat{\Pi} = \prod_v \widehat{G}_v$ we consider the product topology and the Plancherel measure, denoted by $\mu^\Pl$ , given by the product of the local ones. We have so far clarified the settings necessary to properly introduce the measure $\mu$ defined in \eqref{mu}.

\medskip

\begin{rk}
We are interested in the universal family, part of the automorphic dual, henceforth of $\widehat{\Pi}$ which is already endowed with natural topologies. We aim at equidistribution and density results, so we choose among topologies in order to strengthen those properties. We thus seek a quite weak topology, justifying the choice of the product topology instead of the restricted product one, used when discreteness of automorphic forms is sought.
\end{rk}

\subsection{Convergence of $\mu$}
\label{subsec:convergence-mu}
\label{sec:integralscv}

Now that every measure is properly introduced, we come back to the convergence of the integral \eqref{mu}. In order to prove it, it is sufficient to prove the convergence of local integrals. Let us first consider archimedean places. We are able to estimate the integral for we know the involved Plancherel measures \cite{lang_sl2_1985}. The principal series representations with parameter $ir$ have conductor $1+r^2$. Their Plancherel measures are up to a constant $r\mathrm{tanh}(\pi r/2)\dd r$ or $r\mathrm{cotanh}(\pi r/2)\dd r$ according to the parity. The discrete series representation of parameter $k$ has conductor $1+k^2$ and Plancherel measure $k-1$. Hence in all of the three cases, the local integrals converge as do the quantities
\begin{equation}
\int_{0}^\infty \frac{\mathrm{tanh}(\pi r / 2)}{r^3}, \quad \int_{0}^\infty \frac{\mathrm{cotanh}(\pi r / 2)}{r^3} \quad \text{and} \quad \sum_{k \geqslant 1} \frac{k-1}{k^4}.
\end{equation}

As for the finite places where $B$ splits, Section \ref{subsec:constant} computes the associated local integrals which have finite values. The regularization of the integral \eqref{mu} is specifically chosen in order to make the infinite product of those values convergent. The integrals at ramified places are reduced to treating the previous case by the Jacquet-Langlands correspondence, so also converge.

\subsection{Elements of equidistribution}
\label{subsec:elements-equid}
\label{sec:sieve}

Let $S$ be a finite set of places of $B$. Define $F(\widehat{G}_S)$ to be the space of complex bounded functions on $\widehat{G}_S$ supported on a finite number of Bernstein components and whose restriction to the tempered spectrum is continuous outside a set of measure zero for the Plancherel measure restricted to each Bernstein component. Introduce the distribution measure of the truncated universal family,
\[
\mu_Q = \frac{1}{N(Q)} \sum_{\pi \in \mathcal{A}(Q)} \delta_\p, \qquad Q \geqslant 1.
\]

For a measure $\nu$, let 
\begin{equation}
\nu(f) = \int_{\widehat{\Pi}} f(\pi) \dd\nu(\pi) , \qquad f \in F(\widehat{G}_S),
\label{def:equid}
\end{equation}

\noindent  We say that a sequence $(\nu_n)_n$ of Radon positive measures on $\widehat{\Pi}$ weakly converges to a measure $\nu$ if $\nu_n(f)$ converges to $\nu(f)$ for every $f \in F(\widehat{G}_S)$ when $n$ goes to infinity, for every finite set of places $S$.  Since the characteristic functions of relatively quasi-compact open sets of $\widehat{\Pi}$ with zero-measure boundary lie in $F(\widehat{G}_S)$ by the results of  \cite{sauvageot_principe_1997}, this proves that weak convergence of $\mu_Q$ to $\mu$ implies Theorem \ref{thm-equid}. From now on we deal with the measure
\begin{equation}
\nu_Q = \frac{1}{Q^2} \sum_{\pi \in \mathcal{A}(Q)} \delta_\pi, \qquad Q \geqslant 1,
\end{equation}

\smallskip

\noindent easier to handle than $\mu_Q$. This is motivated by the fact, from Theorem \ref{thm-count}, that $N(Q)$ is of asymptotical order $CQ^2$, so that Theorem \ref{thm-equid} is equivalent to: $\nu_Q$ weakly converges to the measure
\begin{equation}
\nu = C\frac{\mu}{\|\mu\|} = \frac{1}{2}\vol(G(F) \backslash G(\A)) \mu.
\end{equation}

\subsection{The Sauvageot density theorem}
\label{subsec:sauvageot}

In order to prove the convergence of $\nu_Q(f)$ to $\nu(f)$ for every function $f \in F(\widehat{G}_S)$, it is sufficient to prove it for Fourier transforms of functions in the Hecke algebra of $G_S$. Indeed, the Sauvageot density theorem \cite{sauvageot_principe_1997} states that any function in $F(\widehat{G}_S)$ can be approximated in that way.
\begin{theorem}[Sauvageot]
Let $S$ be a finite set of places. For every $f \in F(\widehat{G}_S)$ and $\varepsilon > 0$, there exist functions $\phi, \psi \in \mathcal{H}(G_S)$ such that
\begin{itemize}
\item[(i)] $\forall \pi \in \widehat{G}_S, \ |f(\pi) - \widehat{\phi}(\pi)| \leqslant \widehat{\psi}(\pi)$
\item[(ii)] $\mu_S^{\Pl}(\widehat{\psi}) \leqslant \varepsilon$
\end{itemize}
\label{sauvageot}
\end{theorem}

Let us explain how the Sauvageot theorem allows restricting the proof of Theorem \ref{thm-equid} only to functions that are Fourier transforms of functions in the Hecke algebra. Let $f \in F(\widehat{G}_S)$. For $\varepsilon > 0$, there exist $\phi, \psi \in \mathcal{H}(G_S)$ such that $\widehat{\phi}$ and $\widehat{\psi}$ verify the conclusions of the Sauvageot theorem. We then get
\begin{align*}
|\nu_Q(f) - \nu(f)| & \leqslant |\nu_Q(f)-\nu_Q({\widehat{\phi}})| +  |\nu_Q(\widehat{\phi}\,)-\nu({\widehat{\phi}})| +  |\nu(\widehat{\phi}\,)-\nu(f)| \\
& \leqslant \nu_Q(\widehat{\psi}) + |\nu_Q(\widehat{\phi}\,) - \nu(\widehat{\phi}\,)| + \nu(\widehat{\psi}) \\
& \leqslant |\nu_Q(\widehat{\psi}) - \nu(\widehat{\psi})| + 2\nu(\widehat{\psi}) + |\nu_Q(\widehat{\phi}\,) - \nu(\widehat{\phi}\,)|
\end{align*}

\noindent From the definition of $\nu$ and the domination in the Sauvageot theorem it follows, since conductors are at least one, that
\begin{align*}
\nu(\widehat{\psi}) & \ll \zeta_F^{\star}(1) \prod_{v}  \zeta_v(1)^{-1}  \int_{{\widehat{G_v}}} \widehat{\psi}(\pi_v) \frac{\dd\pi_v}{\c(\pi_v)^2}  \\
& \ll  \prod_{v}  \zeta_v(1)^{-1}  \int_{{\widehat{G_v}}} \widehat{\psi}(\pi_v) \dd\pi_v \ll \mu^\Pl_S(\widehat{\psi}) \leqslant \varepsilon
\end{align*}

\noindent So that we get
\begin{equation}
|\nu_Q(f) - \nu(f)| \ll \varepsilon + |\nu_Q(\widehat{\psi}) - \nu(\widehat{\psi})| + |\nu_Q(\widehat{\phi}\,) - \nu(\widehat{\phi}\,)|.
\end{equation}

In order to prove that $\nu_Q$ weakly converges to $\nu$, it is then sufficient to show that the second and third terms vanish for $Q \to \infty$, \textit{i.e.} to prove the theorem for the narrower class of functions $\widehat{\phi}$ and $\widehat{\psi}$. We prove indeed slightly better than what is needed for Theorem \ref{thm-equid}, with a precise asymptotic development in the case of Fourier transforms.

\begin{theorem}
\label{thm:equid-prec}
For every finite set of places $S$ and $\phi \in \mathcal{H}(\widehat{G}_S)$, and every $\varepsilon > 0$,
\begin{equation}
\label{equid-prec}
\nu_Q\left(\widehat{\phi}\,\right) = \nu\left(\widehat{\phi}\,\right) + 
\left\{
\begin{array}{cl}
O\left(Q^{-1+\varepsilon}\right) & \text{ if } F = \Q \text{ and $B$ totally definite;}  \\
O\left(Q^{-\delta_F} \right) & \text{ if $F \neq \Q$} \text{ and $B$ totally definite} ;\\
O\left(\frac{1}{\log Q} \right) & \text{ if $B$ is not totally definite.} \\
\end{array}
\right.
\end{equation}
\end{theorem}

%
%
%

\subsection{Admissible functions}
\label{subsec:selection-S}

Let $\phi \in \mathcal{H} (G_S)$. In this section, we prove that the action of $\widehat{\phi}$ on $\pi$, that is to say $\widehat{\phi}(\pi_S)$, can be assumed to have a selecting effect on the spectral data. The first information is provided by the trace Paley-Wiener theorem of  \cite{bernstein_trace_1986}, that provides the fundamental properties of the Fourier transforms.
\begin{theorem}[Trace Paley-Wiener]
The functions on $\widehat{G}_\p$ lying in the image of the Hecke algebra $\mathcal{H}(G_\p)$ by the Fourier transform are the functions $\widehat{\phi}$ on $\widehat{G}_\p$ such that
\begin{itemize}
\item[(i)] for every standard Levi subgroup $M$ of $G_\p$ and every irreducible representation $\sigma$ of $M$, the function $\psi \mapsto \widehat{\phi}(\mathrm{ind}_{G_\p}^M(\psi \sigma))$ is a regular function on the complex algebraic variety $\psi(M)$ composed of the unramified characters of $M$;
\item[(ii)] there exists an open compact subgroup $K$ of $G_\p$ dominating $\phi$, \textit{i.e.} such that $\phi$ is nonzero only on representations having non trivial $K$-fixed space $\pi^K$.
\end{itemize}
\end{theorem}

It follows that Fourier transforms selects representations of bounded conductor.
\begin{proposition}
\label{prop:S-part-bounds-cond}
Let $S$ be a finite set of finite places, and let $\phi \in \mathcal{H}(G_S)$. There exists $c_\phi > 0$ such that for every generic $\pi \in \widehat{G}_S$ in the support of $\hat{\phi}$, the conductor of $\pi$ is less than $c_\phi$.
\end{proposition}

\proof Since $S$ is a finite set of places it is sufficient to prove the result for a local component. In the case of a finite place $\p$, let $\phi_\p$ be the $\p$-component of $\phi$, where $\p \in S$. The property (ii) of the Trace Paley-Wiener theorem states that its Fourier transform $\widehat{\phi}_\p$  is dominated by a certain open compact subgroup $K$ of $G_\p$, that is to say is supported on representations $\pi_\p$ having nontrivial fixed space $\pi_\p^K$. Since $K$ is open and the sequence $(\overline{K}_{0, \p}(\p^i))_i$ is a filtration in $G_\p$, $K$ contains a certain $\overline{K}_{0, \p}(\p^r)$ and hence $\widehat{\phi}_\p$ is nonzero only for representations of conductor dividing $\p^r$. Since $S$ contains only a finite number of places, this proves that $\widehat{\phi}$ selects only representations $\pi_S \in \widehat{G}_S$ with conductor dividing the product of the corresponding $\p^r$. \qed

\medskip

In order to select automorphic representations in the universal family through trace formula methods, it is necessary to restrict the Fourier transforms considered to the generic spectrum, for otherwise there is no notion of conductor attached to a representation. The following proposition states that it is possible, up to another approximation by density.

\begin{proposition}
\label{PP}
Let $f \in F(\widehat{G}_S)$. Let $\tilde{f}$ be the restriction of $f$ to the generic spectrum, extended by zero elsewhere on the unitary dual. Then $\tilde{f}$ lies in $F(\widehat{G}_S)$.
\end{proposition}

\proof  Recall that the Sauvageot density theorem provides a criterion for functions to be approximated by Fourier transforms. All the properties of the Sauvageot class $F(\widehat{G}_S)$ obviously hold for $\tilde{f}$ except possibly the condition on the discontinuity points. In the case of $\GL(n)$, if $\pi$ is unitary and tempered, then it is generic \cite{prasad_representation_2000}. Therefore the restriction of $f$ and of $\tilde{f}$ to the tempered spectrum coincide, so that if $f$ is continuous on the unitary tempered spectrum thus so is $\tilde{f}$. \qed

\medskip

The conductor of a representation in the generic spectrum is well-defined. Bounding the conductor is not enough for the purposes of the trace formula, it is also necessary to restrict the functions to fixed conductors. However, it is far from obvious that such modified functions are still approximated by Fourier transforms. This is the meaning of the next proposition.

\begin{proposition}
\label{QQ}
Let $S$ be a finite set of finite places and $\q$ an integer ideal supported in $S$. Let $f \in F(\widehat{G}_S)$ supported in the generic spectrum. Let $\bar{f}$ be the restriction of $f$ to the representations of fixed conductor $\q$, extended by zero elsewhere. Then $\bar{f}$ lies in $F(\widehat{G}_S)$.
\end{proposition}
%

\proof We can assume that $f$ is supported on a fixed Bernstein component, say associated with a discrete series representation $\sigma$ on a Levi $M$. We can assume that $\sigma$ is of conductor $\q$, for otherwise the restriction is zero, which clearly lies in $F(\widehat{G}_S)$. The irreducible unitary tempered representations lying in this Bernstein component are given by the fully induced representations $\mathrm{ind}_M^G(\sigma \otimes \chi)$ where $\chi$ is an unramified unitary character. All these induced representations have same conductor since $\chi$ is unitary. Thus the restrictions of $f$ and $\tilde{f}$ to the unitary tempered spectrum coincide. \qed

\medskip

\begin{corollary}
If, for every finite set of places $S$ and every ideal $\q$ supported in $S$, \eqref{equid-prec} holds for functions in $\mathcal{H}(G_S)$ supported in the generic spectrum of fixed conductor $\q$, then \eqref{equid-prec} holds for every function in $\mathcal{H}(G_S)$.
\end{corollary}

\proof Let $\phi \in \mathcal{H}(G_S)$. Let $\bar{\phi}$ be the restriction-extension of $\widehat{\phi}$ to the generic spectrum of fixed conductor $\q$. Propositions \ref{PP} and \ref{QQ} ensure that $\bar{\phi}$ remains in $F(\widehat{G}_S)$. The very same approximation by the Sauvageot density theorem than in the previous section shows that proving the theorem for such functions implies the theorem for every $\widehat{\phi} \in F(\widehat{G}_S)$. \qed

\section{Spectral data}
\label{sec:combi}
\label{subsec:sieve}

\subsection{The Langlands classification}
\label{subsec:Langlands}

The local Langlands classification\index{Langlands classification} of the archimedean admissible dual \cite{knapp_local_1994} of $\GL(2)$ provides a recipe for constructing the admissible representations of reductive groups over archimedian local fields in terms of representations
of its Levi subgroups. Since unitary representations are in particular admissible, it induces a parametrization of the unitary dual of $\GL(2, F_\infty^R)$.  Let $\mathcal{L}_\infty$ the finite set of Levi subgroups of $\GL(2, F_\infty^R)$ containing the diagonal torus. For such a Levi $M$, define $\mathcal{E}_2(M^1)$ to be the set of isomorphism classes of square integrable representations of $M^1$. The only nonempty cases are, at a given archimedean place, 
\begin{itemize}
\item $\mathcal{E}_2(\GL(1, \R)^1)$ consisting in the trivial character and the sign character;
\item $\mathcal{E}_2(\GL(1, \C)^1)$ composed by the characters $z \mapsto z^k/|z|^k$ for integers $k$;
\item $\mathcal{E}_2(\GL(2, \R)^1)$ that is the set of discrete series of weight $k \geqslant 2$.
\end{itemize}

Introduce the set\label{D}\label{delta} $\mathcal{D}$ of $G^R_\infty$-classes of conjugation of pairs $\underline{\delta} = (M, \delta)$ with $M \in \mathcal{L}_\infty$ and $\delta \in \mathcal{E}_2(M^1)$. Write $\h_{M, \C}^\star$ for the trace-zero hyperplane of the complexified dual of the Lie algebra of $M$, which is a finite-dimensional $\C$-vector space. The spectral data $\underline{\delta} \in \mathcal{D}$, consisting of a Levi $M$ in $\mathcal{L}_\infty$ and a discrete series representation $\delta$ in $\mathcal{E}_2(M^1)$, along with $\nu \in \h_{M, \C}^\star$, give rise to an admissible representation of $G_\infty^{R}$ in the following way. The unitary induction $\mathrm{Ind}_P^G(\delta \otimes e^{\nu})$ from the unique parabolic subgroup $P$ containing $M$ is not necessarily irreducible, yet the following holds.
\begin{proposition}[Archimedean Langlands classification]
Let $\underline{\delta} \in \mathcal{D}$ and $\nu \in \h_{M, \C}^\star$, denote $W_\delta$ the stabilizer of $\delta$ in the Weyl group\label{W} of $\h_M$. There is a unique $\nu'$ in the class of $\nu$ modulo translation by $W_\delta$ such that the induction $\mathrm{Ind}_P^G(\delta \otimes e^{\nu'})$ admits an irreducible quotient, that is then unique and denoted by $\pi_{\delta, \nu}$. Moreover, every admissible irreducible representation of $\GL(2, F_\infty^R)$ arises uniquely in this way, up to infinitesimal equivalence.
\end{proposition}

This construction exhausts the admissible dual of $\GL(2, F_\infty^R)$ up to infinitesimal equivalence. This is the archimedean Langlands classification \cite[Theorem 8.54]{knapp_representation_1986}, that can be reformulated as
\begin{equation}
\widehat{\GL}(2)_\infty^{R,1} \simeq \bigsqcup_{\underline{\delta} = (\delta, M) \in \mathcal{D}} \mathcal{E}_2(M^1) \times \mathfrak{h}^\star_{M, \C} / W_\delta,
\end{equation}

\noindent where $\widehat{G}_\infty^{R,1}$ stands for the admissible dual of $\GL(2, F_\infty^R)$ up to infinitesimal equivalence. Note that $\mathcal{D}$ is a discrete set, leading to refer to $\underline{\delta}$ as the discrete archimedean spectral parameter of $\pi$, while $\nu$ in $\h_{M, \C}^{\star}/W_\delta$ is called the continuous archimedean parameter of $\pi$. 

\subsection{Sieving the universal family}

In order to address the problem of the weak convergence of $\nu_Q$ to prove Theorem \ref{thm:equid-prec}, it is necessary to decompose the universal family into smaller sets with fixed spectral data, amenable to trace formula methods. Let $S$ be a finite set of places and  $\phi \in \mathcal{H}(G_S)$. The conductor of $\pi \in \mathcal{A}(G)$ splits into local conductors, in particular can be written
\begin{equation}
\label{conductor-splitting}
c(\pi) = c(\pi_R)  c(\pi^{R}_\infty) c(\pi^R_{S,f})N\c(\pi^{R, S}_f).
\end{equation}

This decomposition emphasizes the different kind of information and behavior each type of place is endowed with, and turns to be a guide for the method. Concerning the split archimedean places, introduce the truncated archimedean split dual
\begin{equation}
\label{Omega}
\Omega(X) = \left\{ \pi^R_\infty \in \widehat{G}^{R}_\infty \ : \ c(\pi^R_\infty) \leqslant X \right\}, \qquad X > 0.
\end{equation}

This set of archimedean parameters factorizes further through the precise Langlands classification recalled in Section \ref{subsec:Langlands}, by fixing discrete spectral parameters, so that
\begin{equation}
\label{omega-splitting}
\Omega(X) = \Omega_{\mathrm{comp}}(X) \sqcup \Omega_{\mathrm{temp}}(X) =  \Omega_{\mathrm{comp}}(X) \sqcup \bigsqcup_{\substack{\underline{\delta} \in \mathcal{D} \\ \underline{\delta} = (M, \delta)}} \Omega_{\underline{\delta}}(X),
\end{equation}

\noindent where
\begin{align*}
\Omega_{\underline{\delta}}(X) & = \left\{ \pi^R_\infty \in \widehat{G}^{R}_\infty \ : \ \exists \nu \in i\h_M^\star, \ \pi^R_\infty \simeq \pi_{\delta, \nu}, \ c(\pi^R_\infty) \leqslant X \right\}\\
\Omega_{\mathrm{comp}}(X) & = \left\{ \pi^R_\infty \in \widehat{G}^{R}_\infty \ : \ \exists \nu \in \h_{M, \C}^\star \backslash i \h_M^\star, \ \pi^R_\infty \simeq \pi_{\star, \nu}, \ c(\pi^R_\infty) \leqslant X \right\} \\
\Omega_{\mathrm{temp}}(X) & = \bigsqcup_{\substack{\underline{\delta} \in \mathcal{D} \\ \underline{\delta} = (M, \delta)}} \Omega_{\underline{\delta}}(X) 
\end{align*}

\noindent and the notation $\simeq \pi_{\star, \nu}$ stands for the existence of a $\underline{\delta} \in \mathcal{D}$ such that the representation is isomorphic to $\pi_{\delta, \nu}$. The set $\Omega_\mathrm{comp}$ is called the complementary part of the archimedean spectrum, while the remaining part is the tempered part of the spectrum. This denomination is motivated by the fact that the representation $\pi_{\delta, \nu}$ is tempered if and only if $\nu$ lies in $i\h_M^\star$.

As for the remaining places, given a finite set of places $S$, recall that every ideal $\m$ is decomposed in the form $\m = \m_S\m^S$, where such a decomposition always means that $\m^S$ is the prime-to-$S$ part of $\m$, \textit{i.e.} is such that $\m^S \wedge S = 1$, and $\m_S$ if the $S$-part of $\m$, \textit{i.e.} satisfies $\mathrm{supp}(\m_S) \subseteq S$. The same decomposition is used without further notice for the other letters. The multiplicative conductor of the finite split places is fixed to a certain ideal $\q$ of $\mathcal{O}^R$, and the isomorphism class of the ramified part is fixed to a certain isomorphism class $\sigma_R \in \widehat{G}_R$.

Recall from Proposition \ref{QQ} that the function $\bar{\phi}$ fixes the conductor of the finite $S$-component to be equal to a certain $\q_S$. Thus, according to \eqref{conductor-splitting} and the choices made above, the universal family $\mathcal{A}(Q)$ decomposes as
\begin{equation}
\label{uf-decomposition}
\bigsqcup_{\substack{N\q \leqslant Q\\ \q \wedge R = 1}} \bigsqcup_{\substack{\sigma_R \in \widehat{G}_R \\ c(\sigma_R) \leqslant Q/N\q}} \bigsqcup_{\substack{\underline{\delta}\in\mathcal{D} \\ \underline{\delta} = (M, \delta)}} \mathcal{A}(\q, \sigma_R, \delta, Q) \sqcup \bigsqcup_{\substack{N\q \leqslant Q\\ \q \wedge R = 1}} \bigsqcup_{\substack{\sigma_R \in \widehat{G}_R \\ c(\sigma_R) \leqslant Q/N\q}} \mathcal{A}_{\mathrm{comp}}(\q, \sigma_R, Q) ,
\end{equation}

\noindent where the sets of fixed spectral data are
\begin{align*}
\mathcal{A}(\q, \sigma_R, \delta, Q) & =  \left\{ \pi \in \mathcal{A}(G) \ : \ \pi_R \simeq \sigma_R, \ \c(\pi^{R}_f) = \q, \ \pi_\infty^R \in \Omega_{\underline{\delta}}\left(\frac{Q}{N\q c(\sigma_R)}\right) \right\} \\
\mathcal{A}_{\mathrm{comp}}(\q, \sigma_R, Q) & =  \left\{ \pi \in \mathcal{A}(G) \ : \ \pi_R \simeq \sigma_R, \ \c(\pi^{R}_f) = \q, \ \pi_\infty^R \in \Omega_{\mathrm{comp}}\left(\frac{Q}{N\q c(\sigma_R)}\right) \right\}
\end{align*}

\noindent This decomposition \eqref{uf-decomposition} of the universal family reduces the study of the whole family to the harmonic families $\mathcal{A}(\q, \sigma_R, \delta, Q)$, easier to grasp in the context of trace formulas. What is critical is to having got rid of the condition of belonging to $\mathcal{A}(Q)$, decomposed in local conditions. It induces a decomposition of the counting measure as
\begin{equation}
\label{decomposition1}
\begin{split}
\nu_Q(\widehat{\phi}\,) & =  \frac{1}{Q^2} \sum_{\pi \in \mathcal{A}(Q)} \widehat{\phi}(\pi)  \\
& =  \frac{1}{Q^2} \sum_{\substack{\pi \in \mathcal{A}(G) \\ c(\pi_R)c(\pi^{R, S}_f)c(\pi^R_{S,f})c(\pi^R_\infty) \leqslant Q}} \widehat{\phi}(\pi)  \\
& = \frac{1}{Q^2}  \sum_{\substack{\sigma_R \in \widehat{G}_R \\ c({\sigma_R}) \leqslant Q}}  \sum_{\substack{N\q \leqslant Q/c(\sigma_R) \\ \q \wedge  R = 1}} \sum_{\substack{\underline{\delta}\in\mathcal{D} \\ \underline{\delta} = (M, \delta)}}
\sum_{\substack{\pi \in \mathcal{A}(\q, \sigma_R, \delta, Q) }} \widehat{\phi}(\pi) \\
& \qquad \qquad \qquad + \frac{1}{Q^2} \sum_{\substack{\sigma_R \in \widehat{G}_R \\ c({\sigma_R}) \leqslant Q}}  \sum_{\substack{N\q \leqslant Q/c(\sigma_R)\\ \q \wedge  R = 1}} 
\sum_{\substack{\pi \in \mathcal{A}_{\mathrm{comp}}(\q, \sigma_R, Q)}} \widehat{\phi}(\pi) 
\end{split}
\end{equation}

\noindent where the sum over $\q$ is meant to run through ideals of $\mathcal{O}^{R}$. The complementary part corresponds to the second sum appearing in the line above and will be dealt with later and shown to contribute as an error term. Denote ${A}(\q, \sigma_R, \delta, Q;\PPP)$ and ${A}_{\mathrm{comp}}(\q, \sigma_R, Q;\PPP)$ the innermost parts of the splitting in the first summation above, that is to say
\begin{equation}
\label{A}
A(\q, \sigma_R, \delta, Q;\PPP) = \sum_{\substack{\pi \in \mathcal{A}(\q, \sigma_R, \delta, Q)}} \widehat{\phi}(\pi), 
\end{equation}

\noindent and analogously for $A_{\mathrm{comp}}(\q, \sigma_R, Q;\PPP)$.

\subsection{Old and new forms}
\label{oldnew}

The universal family \eqref{UF} sees no multiplicities\index{multiplicities}, but the trace formula counts them. The spectral multiplicities associated to the spectral decomposition of $L^2(G(F) \backslash G(\A))$, which are more suitable weights for the forthcoming computations, are given by 
\begin{equation}
\label{multiplicities}
m\left(\pi, \q\right) = \dim \left(\pi^{\overline{K}_0(\q)} \right),
\end{equation}

\noindent where
\begin{equation}
Z K_0(\q) = \prod_{\substack{\p^r || \q}} Z_\p K_{0, \p}\left(\p^r\right) \subseteq B^{\times}\left(\A^R_f\right),
\end{equation}

\noindent and $\overline{K}_0(\q)$ stands for the image of $ZK_0(\q)$ under the natural projection $B^\times \rightarrow G$. The choice is made so that $m(\pi, \q) \neq 0$ is equivalent to $\c(\pi^R_f) \div \q$. The analogous sum to \eqref{A} additionally weighted by the multiplicities is
\begin{equation}
\label{B}
B(\q, \sigma_R, \delta, Q;\PPP) = \sum_{\substack{\pi \in \mathcal{B}(\q, \sigma_R, \delta, Q)}} m\left(\pi^{S}, \q^S\right) \widehat{\phi}(\pi), 
\end{equation}

\noindent where
\begin{equation*}
\mathcal{B}(\q, \sigma_R, \delta, Q) = \left\{ \pi \in \mathcal{A}(Q) \ : \ \pi_R \simeq \sigma_R, \  \c\left(\pi^{R}_f\right) \div \q, \ \pi^R_\infty \in \Omega_{\underline{\delta}}(Q/N\q c(\sigma_R)) \right\}.
\end{equation*}

\medskip

\noindent and analogously for $B_{\mathrm{comp}}(\q, \sigma_R, Q;\PPP)$. The sum defined by \eqref{A} counts the newforms\index{newforms} while \eqref{B} counts the old ones with respect to finite prime-to-$S$ split places. The relation between them lies in the following lemma.

\begin{lemma}  
\label{lem:sieve}
Let $\q$ prime to $R$,  $\sigma_R$ an irreducible unitary representation of $G_R$, $\underline{\delta} \in \mathcal{D}$ and $\phi \in \mathcal{H}(G_S)$. Let $\lambda_2= \mu \star \mu$ where $\mu$ is the Möbius function. For every $Q \geqslant 1$,
\begin{align*}
A\left(\q,\sigma_R, \delta, Q; \PPP\right) & = \sum_{\d\div \q} \lambda_2\left(\frac{\q}{\d}\right) B\left(\d, \sigma_R, \delta, Q; \PPP\right) \\
A_{\mathrm{comp}}\left(\q,\sigma_R, Q; \PPP\right) & = \sum_{\d\div \q} \lambda_2\left(\frac{\q}{\d}\right) B_{\mathrm{comp}}\left(\d, \sigma_R, Q; \PPP\right) \\
\label{sieve}
\end{align*}
\end{lemma}

\proof Recall that, for every finite split place $\p$, Casselman gives the local multiplicites
\begin{equation}
\dim \sigma_{\p}^{K_0\left(\p^{\f(\sigma_\p)+i}\right)} = i+1, \qquad i \geqslant 0.
\end{equation}

\noindent From this immediately follows, after taking the product over all finite split places, that the global multiplicities are
\begin{equation}
m\left(\sigma, \q\right) = \tau_2\left(\frac{\q}{\c(\sigma^{R}_f)}\right),
\label{mult}
\end{equation}

\noindent where $\tau_2 = 1 \star 1$ is the divisor function. Since $\left(\sigma^{R}\right)^{\overline{K}_0\left(\q\right)} \neq 0$ implies $\c(\sigma^{R}) \div \q$, the sum defining $B(\q, \sigma_R, \delta, Q;\PPP)$ is eventually reduced to a sum over $\c(\sigma^{R}) \div \q$. Thus, by the precise knowledge \eqref{mult} of the multiplicities, 
\begin{equation}
\begin{split}
B\left(\q,\sigma_R, \delta, Q; \PPP\right) &  = \sum_{\d\div \q} \sum_{\sigma \in \mathcal{A}(\d, \sigma_R,\delta, Q)} \tau_2\left(\frac{\q}{\c(\sigma^{R}_f)}\right) \widehat{\phi}(\sigma)  \\
&  = \sum_{\d\div \q} \tau_2\left(\frac{\q}{\d}\right) \sum_{\sigma \in \mathcal{A}(\d, \sigma_R, \delta, Q)} \widehat{\phi}(\sigma) \\
&  = \sum_{\d\div \q} \tau_2\left(\frac{\q}{\d}\right) A\left(\d,\sigma_R, \delta, Q; \PPP\right)
\end{split}
\end{equation}

\noindent so that $B = \tau_2 \star A$, with a slight abuse of notation. Hence, by Möbius inversion,
\begin{equation}
A\left(\q,\sigma_R, \delta, Q; \PPP\right) = \sum_{\d\div \q} \lambda_2\left(\frac{\q}{\d}\right) B\left(\d,\sigma_R, \delta, Q; \PPP\right),
\end{equation}

\noindent achieving the first part of the claim. The proof carries on, \textit{mutatis mutandis}, for the quantities relative to the complementary spectrum. \qed

\medskip

Summing over the spectral data appearing in the decomposition \eqref{decomposition1}, the counting measure rewrites as
\begin{equation}
\label{step3}
\nu_Q(\widehat{\phi}\,) = \nu_{\mathrm{temp}, Q}(\widehat{\phi}\,) + \nu_{\mathrm{comp}, Q}(\widehat{\phi}\,), \\
\end{equation}

\noindent where
\begin{align*}
\nu_{Q, \mathrm{temp}}(\widehat{\phi}\,) & =  \frac{1}{Q^2}  \sum_{\substack{\sigma_R \in \widehat{G}_R \\ c({\sigma_R}) \leqslant Q}} \sum_{\substack{N\q \leqslant Q/c(\sigma_R)) \\ \q \wedge R = 1}}   \sum_{\substack{\underline{\delta} \in \mathcal{D} \\ \underline{\delta} = (M, \delta)}} \sum_{\d \div \q} \lambda_2 \left(\frac{\q}{\d}\right) B\left(\d, \sigma_R, \delta, Q; \PPP\right) \\
\nu_{Q, \mathrm{comp}}(\widehat{\phi}\,) & =  \frac{1}{Q^2}  \sum_{\substack{\sigma_R \in \widehat{G}_R \\ c({\sigma_R}) \leqslant Q}} \sum_{\substack{N\q \leqslant Q/c(\sigma_R)) \\ \q \wedge R = 1}}    \sum_{\d \div \q} \lambda_2 \left(\frac{\q}{\d}\right) B_{\mathrm{comp}}\left(\d, \sigma_R, Q; \PPP\right) \\
\end{align*}

\section{Trace formula}
\label{sec:fts}
\label{sec:trace-formula}

Trace formul\ae{} give relations between spectral and geometrical quantities, the latter being often easier to manipulate. We present here the Selberg trace formula and express the sought tempered old forms numbers $B(\d, \sigma_R, \delta, Q; \PPP)$ and its complementary counterpart $B_{\mathrm{comp}}(\d, \sigma_R, Q; \PPP)$ as a spectral side of this trace formula for a suitable test function, leaving us with the geometric side to estimate.

\subsection{Selberg trace formula}
\label{subsec:fts}

Since the automorphic quotient of $G$ is compact, the original formulation of the trace formula\index{trace formula}, due to Selberg in 1956, can be used and combined with the multiplicity one theorem. If $\Phi$ is a function in the Hecke algebra $\mathcal{H}(G(\A))$, then
\begin{equation}
\label{fts}
J_{\geom}(\Phi)= J_{\spec}({\Phi}),
\end{equation}

\noindent where the spectral and geometrical parts\index{trace formula!spectral part}\index{trace formula!geometrical part} are as follows. The geometrical part is

\begin{equation}
\label{fts:geometrical-part}
J_{\geom}(\Phi) := \sum_{\{\gamma\}} \vol\left(G_\gamma(F) \backslash G_\gamma(\A)\right) \int_{G_\gamma(\A) \backslash G(\A)} \Phi\left(x^{-1}\gamma x\right) \dd x .
\end{equation}

\noindent The sum runs through conjugacy classes $\{\gamma\}$ in $G(F)$. Since $\Phi$ is compactly supported and $G(F)$ is discrete, the sum is finite. However its length depends on the support of $\Phi$ what turns to be a critical difficulty for estimations, for this support will depend on the spectral parameters. The inner integrals appearing in this geometric side are called the orbital integrals\index{orbital integral, $\mathcal{O}_\gamma(\Phi)$}, defined by
\begin{equation}
\label{OI}
\mathcal{O}_\gamma(\Phi) = \int_{G_\gamma(\A) \backslash G(\A)} \Phi\left(x^{-1}\gamma x\right) \dd x .
\end{equation}

The spectral part is
\begin{equation}
\label{fts:spectral-part}
  J_{\spec}({\Phi})  = \sum_{\substack{\pi \subseteq L^2(G(F) \backslash G(\A))}} m(\pi) \widehat{\Phi}(\pi).
\end{equation}

\noindent Here $\pi$ go through the isomorphism classes of unitary irreducible subrepresentations of $G(\A)$ in $L^2(G(F) \backslash G(\A))$, and recall that $\widehat{\Phi}$ is the Fourier transform of $\Phi$, see Section \ref{sec:plancherel-formulas}.

\medskip

\begin{rk}
The formulaton of the spectral part \eqref{fts:spectral-part} is Selberg's original one. The weights $m(\pi)$ are the multiplicities of the $\pi$'s in the discrete part of the spectral decomposition of $L^2(G(F) \backslash G(\A))$. The multiplicity one theorem ensures these to be less than one, and the indexation by $\pi$ actually part of $L^2(G(F) \backslash G(\A))$ makes them nonzero, hence equal to one.
\end{rk}

\medskip

The admissible dual can be decomposed into tempered\label{Jtemp} representations and non-tempered representations. In view of \eqref{omega-splitting} and anticipating that the selecting function at split archimedean places behaves differently on the tempered spectrum and on the complementary one, it is natural to introduce the tempered and complementary spectral parts as
\begin{align*}
  J_{\mathrm{temp}}({\Phi})  & = \sum_{\substack{\pi \subseteq L^2(G(F) \backslash G(\A)) \\ \pi_\infty^R \simeq \pi_{\star, \nu} \\ \nu \in \Omega_{\mathrm{temp}}}} m(\pi) \widehat{\Phi}(\pi) \\
    J_{\mathrm{comp}}({\Phi})  & = \sum_{\substack{\pi \subseteq L^2(G(F) \backslash G(\A)) \\ \pi_\infty^R \simeq \pi_{\star, \nu} \\ \nu \in \Omega_{\mathrm{comp}}}} m(\pi) \widehat{\Phi}(\pi) 
\end{align*}

As announced in the outlook of the method, in order to have a problem amenable to the trace formula it is necessary to formulate statistics quantities on the universal family as a spectral side,  hence needed to select it by the Fourier transforms of suitable test functions. The aim of the present section is to construct a function $\Phi \in \mathcal{H}(G)$ such that, up to an error term,
\begin{align}
\label{selection}
J_{\mathrm{temp}}(\Phi) & = B\left(\d, \sigma_R, \delta, Q; \PPP\right) \\
J_{\mathrm{comp}}(\Phi) & = B_{\mathrm{comp}}\left(\d, \sigma_R, Q; \PPP\right) 
\end{align}

In the case of factorizable test functions $\Phi = \otimes_v \Phi_v$, the spectral side of the trace formula factorizes as
\begin{equation}
\label{Fourier-transform:factorization}
\widehat{\Phi}(\pi) = \prod_v \widehat{\Phi}_v (\pi_v).
\end{equation}

Hence, in order to achieve the spectral selection \eqref{selection} it is sufficient locally select the conditions appearing in the decomposition of the universal family \eqref{step3} through Fourier transforms. The following sections are dedicated to construct local test functions doing so, aim reached in Lemma \ref{lem:selecting}. The places of $F$ fall into four categories:
\begin{itemize}
\itex the split finite part, corresponding to $\p \notin R \cup S$, where the arithmetic conductor is caught by the means of an explicit filtration, see Section \ref{filtration};
\itex the split finite part in the support of the test function $\widehat{\phi}$, corresponding to $\p \in S \backslash R$, whose conductor is fixed by $\widehat{\phi}$, see Proposition \ref{QQ};
\itex the ramified part, corresponding to the finite number of $v \in R$, which is handled by fixing the representations at those places by means of matrix coefficients;
\itex the split archimedean part, parametrized by spectral data that are handled by selecting functions provided by Paley-Wiener theorems.
\end{itemize}

\subsection{Selecting the split conductor}

For an ideal $\d$ of $\mathcal{O}$, introduce the congruence subgroup given by the product of the corresponding local congruence subgroups in \eqref{filtration}, that is to say
\begin{equation}
\label{filtration-global}
K_0(\d) = \prod_{\p^r || \d} K_{0, \p}(\p^r).
\end{equation}

The following result gives a test function whose Fourier transform selects the finite split conductor.

\begin{lemma}
\label{lem:ft-finite-split-conductor}
For an ideal $\d$ of $\mathcal{O}$, let
\begin{equation}
\label{ft-finite-split-conductor}
\varepsilon_{\d} = \vol\left(\overline{K}_0(\d)\right)^{-1} \mathbf{1}_{\overline{K}_0(\d)}.
\end{equation}
Its Fourier transform selects the multiplicity relative to $\d$. More precisely, 
\begin{equation}
\label{Fourier-ft-finite-split-conductor}
\widehat{\varepsilon}_\d (\pi) = m(\pi, \d), \qquad \pi \in \mathcal{A}(G).
\end{equation}
\end{lemma}

\proof Let $\pi$ be an automorphic representation of $G$. Then $\pi(\varepsilon_{\d})$ is the projection of the representation space $V_{\pi}$ on the subspace $\pi^{\d}$ of the fixed vectors by $\overline{K}_0(\d)$ under the action of $\pi$. Indeed, every $\pi(\varepsilon_{\d})v$, for $v$ in  $V_{\pi}$, is $\overline{K}_0(\d)$-invariant, for it is an averaging over the action of $\overline{K}_0(\d)$. For $k_0 \in \overline{K}_0(\d)$ and $v \in V_{\pi}$, 
\begin{align*}
  \pi(k_0) \pi\left(\varepsilon_{\d}\right)v & = \vol\left(\overline{K}_0(\d)\right)^{-1} \pi(k_0) \int_{\overline{K}_0(\d)} \pi(k) v \dd k   \\
& = \vol\left(\overline{K}_0(\d)\right)^{-1} \int_{\overline{K}_0(\d)} \pi(k_0 k) v \dd k  \\
& = \vol\left(\overline{K}_0(\d)\right)^{-1} \int_{\overline{K}_0(\d)} \pi(k) v \dd k  = \pi\left(\varepsilon_{\d}\right)v
\end{align*}

\noindent so that its image lies in $\pi^{\d}$. The action of $\pi(\varepsilon_{\d})$ is also idempotent, more precisely the identity on $\pi^{\d}$. Indeed, for $v_0 \in \pi^{\d}$,
\begin{align*}
\pi\left(\varepsilon_{\d}\right) v_0 & =  \vol\left(\overline{K}_0(\d)\right)^{-1} \int_{\overline{K}_0(\d)} \pi(k) {v_0} \dd k  \\
& = \vol\left(\overline{K}_0(\d)\right)^{-1} \int_{\overline{K}_0(\d)} v_0 \dd k \\
& = v_0
\end{align*}

\noindent Hence $\pi(\varepsilon_{\d})$ is an idempotent endomorphism of image $\pi^{\d}$, \textit{i.e.} a projection on $\pi^{\d}$. The trace of a projection is its rank, that is to say $\widehat{\varepsilon}_{\d}(\pi)$ is the dimension of the fixed vector spaces $\pi^{\d}$. Those are the sought multiplicities $m(\pi, \d)$, in particular are nonzero if and only if $\c(\pi) \div \d$. \qed

\subsection{Selecting the ramified part}

For ramified places, less is known concerning the representations and the choice made in the decomposition \eqref{step3} is to fix the corresponding isomorphism class. In the finite dimensional case, knowing matrix coefficients is sufficient to determine the underlying matrix. This property still holds \cite[Corollary 10.26]{knightly_traces_2006} for supercuspidal representations in the following sense. Let $\sigma_R$ be a unitary representation of $G_R$. A matrix coefficient\index{matrix coefficient, $\xi_\sigma$} associated to $\sigma_R$ is a function of the form, given $v$ and $w$ in the space of $\sigma_R$, 
\begin{equation}
\label{matrix-coeff:def}
\begin{array}{cccc}
\displaystyle \xi_{\sigma_R}^{v, w} : &  G_R & \longrightarrow & \C \\
& g & \longmapsto & \langle \sigma(g)v, w \rangle
\end{array}
\end{equation}

Matrix coefficients are continuous functions on $G_R$, are compactly supported since $G_R$ is compact, and are locally constant at finite places and smooth at archimedean places.

\medskip

\rk The fact that matrix coefficients is considered only for ramified places is critical for selecting purposes. The loss of the compactness of the support for matrix coefficients in the split case, where some automorphic representations are not supercuspidal, make them fail to select the corresponding isomorphism class. Such a purpose can be achieved by means of existence theorem, yet are less precise, see \cite{knightly_traces_2006}. This is the reason why the non-totally definite case or the $\GL(2)$ case are analytically harder to deal with, see Section \ref{BM-tf}.

\medskip
\begin{proposition}
Let $\sigma$ and $\pi$ be automorphic representations of $G_R$, and introduce $d_\pi$ the formal degree of $\pi$. Then for every unit vectors $v$ and $w$ in the representation space of $\sigma$, 
\begin{equation}
\label{matrix-coeff:orthog-rel}
\pi\left(\xi_{\sigma}^{v, w}\right) w = \mathbf{1}_{\pi \simeq \sigma} \frac{\langle w, v \rangle }{d_\pi} v.
\end{equation}
\end{proposition}

Taking for $v$ a vector of norm $d_\pi^{1/2}$, it follows that $\pi\left( \xi_\sigma^{v,v}\right)$ is the orthogonal projection onto $\C v$ and in the meanwhile selects the $\pi$'s isomorphic to $\sigma$. Considering its trace, this can be restated as follows.
\begin{proposition}
\label{prop:matrix-coeff}
Let $\sigma$ and $\pi$ be automorphic representations of $G_R$. Let $v$ be a vector of norm one in the representation space of $\sigma$. Then,
\begin{equation}
\label{matrix-coeff:selecting-function}
\widehat{\xi_\sigma^{v,v}}(\pi) = \mathbf{1}_{\pi \simeq \sigma}.
\end{equation}
\end{proposition}

From now on, denote $\xi_\sigma$ any choice of matrix coefficient as in Proposition \ref{prop:matrix-coeff}.

\subsection{Approximate localizing at archimedean split places}
\label{BM-tf}
\label{subsec:arch-params}

Based on the decomposition of the universal family \eqref{uf-decomposition}, for a bounded set of continuous parameters $\Omega$ in $\mathfrak{h}_{M, \C}^\star$, the question is to select representations lying in sets of the form
\begin{align*}
\mathcal{A}(\q, \sigma_R, \delta, Q) & =  \left\{ \pi \in \mathcal{A}(Q) \ : \ \pi_R \simeq \sigma_R, \ \c(\pi^{R}_f) = \q, \ \pi_\infty^R \in \Omega_{\underline{\delta}} \right\}.
\end{align*}

This is a feature of non-compact archimedean groups: their duals are no more discrete and hence admit a continuous parametrization. However, many tools in harmonic analysis on Lie groups rely on narrow classes of functions among which characteristic functions of such sets $\Omega_{\underline{\delta}}$ are not, requiring a smoothing construction to get admissible functions lying near them. This procedure is provided by the fundamental work of \cite{DKV}. \cite{brumley_counting_2018} adapted this method to the automorphic setting on $\GL(n)$ and construct a function localizing around spectral parameters $(\underline{\delta}, \nu)$, where $\underline{\delta}$ is a fixed archimedean discrete spectral datum and where $\nu$ is a continuous parameter lying into a bounded set of tempered parameters $\Omega_{\underline{\delta}}$.  Such smoothing procedures behave well on tempered parameters, leading to assume $\Omega$ to be a bounded set of tempered parameters of fixed discrete part $\underline{\delta}$, leaving the non-tempered part of $\Omega$ to be proven negligible compared to the tempered contribution, see Section \ref{sec:COMP}.

Introduce $\phi$, which in this section is a function in the Hecke algebra of $G^R_\infty$ and should be denoted $\phi_{\infty}^R$ in the following ones. The aim is to find a smooth enough function for trace formula purposes approximating the characteristic function of $\Omega_{\underline{\delta}}$.  \cite{brumley_counting_2018} achieve this goal, constructing a function  $h_\rho^{\delta, Q}$ of Paley-Wiener class with exponential type $\rho > 0$ as tempered spectral-localizing function. Let
\begin{equation}
\label{hdeltaXphi}
h_\rho^{\delta, Q, \phi}:= \widehat{\phi} \ h_\rho^{\delta, Q}, \qquad \delta \in \mathcal{D}.
\end{equation}

A direct consequence of their result is the following lemma, where the remainder term is willingly hidden in order to ease the exposition. What is of critical importance are the bounds on this undisclosed error term, precisely stated in Lemma \ref{lem:BM-bounds}.

\begin{lemma}
\label{lem:fta}
\label{lem:arch-tf-spectral}
For every discrete spectral datum $(M, \delta) \in\mathcal{D}$, there is a function $\epsilon_\rho^{\delta, Q}$ such that for every $(M, \tau) \in \mathcal{D}$ and $\nu \in \mathfrak{h}_{M, \C}^\star$,
\begin{align*}
(i) \qquad & h_\rho^{\delta, Q, \phi} (\tau, \nu) = \mathbf{1}_{\substack{\tau \in W_\delta \\ \nu \in \Omega}}  \ \widehat{\phi}(\tau, \nu) + \epsilon_\rho^{\delta, Q}(\tau, \nu) \\
(ii) \qquad & h_\rho^{\delta, Q, \phi} (\tau, \nu) \ll \mathbf{1}_{\substack{\tau \in W_\delta \\ \Re(\nu) \in \Omega}} e^{\rho \|\Re \nu\|} 
\end{align*}
\end{lemma}

\proof This is just encapsulating the results of \cite[Lemma 9.2]{brumley_counting_2018} and multiplying them by $\widehat{\phi}$.  \qed

\medskip

The same arguments used by \cite{brumley_counting_2018} hold with $h_\rho^{\delta, Q}$ replaced by $h_\rho^{\delta, Q, \phi}$. In particular, a version of the Paley-Wiener proven by  \cite{clozel_theoreme_1990} provides a function $f_\rho^{\delta, Q, \phi}$\label{fdeltaXphi} whose Fourier transform is $h_\rho^{\delta, Q, \phi}$. 

\subsection{The chosen test function}
\label{subsec:fonction-test}

The weighted counting number $B(\d, \sigma_R,  \delta, Q;  \PPP)$ should be written as a spectral side in the trace formula. Introduce the test function
\begin{equation}
\Phi_{\d, \pi_R, \delta, Q, \rho; \PPP} = \prod_v \Phi_v, 
\label{def:ft}
\end{equation}

\noindent which is built with the following local functions:
\begin{center}
\begin{tabular}{|L|L|}
\hline
\text{Places $v$} & \text{Local test function $\Phi_v$}\rule{0pt}{2.6ex} \\
\hline 
\notin S, \notin R, <\infty & \varepsilon_{\d, v}\rule{0pt}{2.6ex} \\
\hline
\notin S, \notin R, \in \infty & f^{\delta, Q}_{\rho, v}\rule{0pt}{2.6ex} \\
\hline
\notin S, \in R & \xi_{\pi_v}\rule{0pt}{2.6ex} \\
\hline
\in S, \notin R, <\infty & \phi_v \rule{0pt}{2.6ex} \\
\hline
\in S, \notin R, \in \infty & f^{\delta, Q, \phi}_{\rho, v} \rule{0pt}{2.6ex} \\
\hline
\in S, \in R &  \xi_{\pi_v}\widehat{\phi_v}(\pi_v)\rule{0pt}{2.6ex} \\
\hline
\end{tabular}
\end{center}

\noindent where
\begin{itemize}
\itex $\phi_v$ is the local component of $\PPP$ on $G_v$;
\itex $\xi_{\pi_v}$ is a matrix coefficient for $\pi_v$;
\itex $\varepsilon_{\d}$ is the function introduced in Lemma \ref{lem:ft-finite-split-conductor}, $\varepsilon_{\d, v}$ its $v$-component;
\itex $f_\rho^{\delta, Q, \phi}$ is the function constructed Lemma \ref{lem:arch-tf-geometric}, with $\Omega = \Omega(Q/N\q c(\pi_R))$.
\end{itemize}

\medskip 

The sought weighted measure is barely reached by the spectral side with $\Phi_{\d, \pi_R, \delta, Q, \rho; \PPP}$, as stated in the following lemma.

\begin{lemma}
\label{lem:selecting} Let $Q \geqslant c_\phi$. Let $\d \wedge R = 1$, $\pi_R \in \widehat{G}_R$, $\delta \in \mathcal{E}_2(M^1)$ for an $M \in \mathcal{L}_\infty$. Then
\begin{equation}
B\left(\d, \pi_R, \delta, Q;  \PPP\right) = J_{\mathrm{temp}}\left({\Phi}_{\d, \pi_R, \delta, Q, \rho; \PPP}\right) + O(\Xi(\PPP, \pi_R)) + O(\partial_\rho B(\d, \pi_R, \delta, Q)),
\label{ftsapplied}
\end{equation}

\noindent where, introducing the set $X^\mathrm{ur}(G)$ of unramified characters of $G(\A)$,
\begin{equation}
 \Xi(\PPP, \pi_R)  = \sum_{\substack{\chi \in X^\mathrm{ur}(G) \\ \chi_R \simeq \pi_R}} m(\chi^{R}, \d) \widehat{\phi}(\chi),
\end{equation}

\noindent and
\begin{equation}
\label{blouh}
\partial_\rho B(\d, \pi_R, \delta, Q)  = \int_{\substack{\pi \in \mathcal{A}(\d, \pi_R, \delta) \\ \nu \in i\mathfrak{h}_M^\star}} \tau_{2}\left( \frac{\d}{\c(\pi^{R})} \right) \epsilon_\rho^{\delta, Q}(\tau, \nu) \dd\nu, 
\end{equation}

\noindent where this last integral means an integration over $\pi \in \mathcal{A}(G)$ of fixed discrete spectral data $\d$, $\pi_R$ and $\delta$, and with continuous parameters varying in $i\mathfrak{h}_M^\star$.
\end{lemma}

\proof Let $\Phi = \Phi_{\d, \pi_R, \delta, Q, \rho; \PPP}$. In order to determine the Fourier transform of $\Phi$ recall that for every places $v$, $w$ and every $a \in \mathcal{H}(G_{v, w})$, $\widehat{a_va_w} = \widehat{a}_v\widehat{a_w}$. Thus,
\begin{equation}
\widehat{\Phi} = \prod_v \widehat{\Phi}_v = h^{\delta, Q, \phi}_\rho \prod_{v \in R} \widehat{\xi}_{\pi_v} \prod_{\substack{\p \notin R \\ \p < \infty \\ \p \notin S \\ \p^r || \d}} \widehat{\varepsilon}_{\p^r, v} \prod_{\substack{\p \notin R \\ \p < \infty \\ \p \in S}} \widehat{\phi}_\p.
\label{transform:ft}
\end{equation}

\medskip
Hence only the Fourier transforms of the local components of the test function have to be determined.  The finite prime-to-$S$ split part $\varepsilon_\d$ is shown to transform into the characteristic function of conductors dividing $\d$ in Lemma \ref{lem:ft-finite-split-conductor} weighted by the corresponding multiplicities. The ramified local parts $\xi_{\pi_v}$ are known to transform into the characteristic functions of the isomorphism class of $\pi_v$ by Lemma \ref{prop:matrix-coeff}. The transform of the archimedean split part is shown to approximate the selecting function of bounded conductors in Lemma \ref{lem:arch-tf-spectral}, up to a smoothing error term $\epsilon_\rho^{\delta, Q}$.  The action of the Fourier transform of $\Phi$ on the tempered part follows, and \eqref{transform:ft} yields, for $\sigma \in \mathcal{A}(G)$ with archimedean split parameters $(\tau,  \nu)$,
\begin{equation}
\label{ft-Fourier-transformed}
\widehat{\Phi}(\sigma) =  m(\sigma^{R}, \d) \widehat{\phi}(\sigma_f)  \mathbf{1}_{\substack{\sigma_R \simeq \pi_R \\ \c(\sigma^{R}) \div \d}}  \left( \mathbf{1}_{\substack{\tau \in W_\delta \\ \nu \in \Omega}}  \ \widehat{\phi}(\tau, \nu) + \epsilon_\rho^{\delta, Q}(\tau, \nu) \right).
\end{equation}

Nevertheless, these conditions also stand for characters: in order to not being killed by $\widehat{\Phi}$ they have to be trivial on $\overline{K}_0(\d)$, \textit{i.e.} they have to be unramified since $\det(\overline{K}_0(\d)) = \mathcal{O}^R$. Moreover, they have to be isomorphic to $\pi_R$ at ramified places. The Fourier transform of the chosen test function hence does not vanish on unramified characters, unlike awaited. The corresponding extra contribution $\Xi$ is treated separately in Lemma \ref{lem:characrers}, for characters are easier to embrace and it will be shown to contribute as an error term.  

After integrating over the tempered spectrum, it follows by roughly bounding $\widehat{\phi}$ in the remainder smoothing term,
\begin{align*}
\label{Jtemp}
J_{\mathrm{temp}}(\Phi) & = \sum_{\substack{\sigma \in \mathcal{B}(\d, \pi_R, \delta, Q)}} \tau_{2}\left( \frac{\d}{\c(\sigma^{R})}\right) \widehat{\phi}(\sigma) \\
& \qquad + O\left( \int_{\substack{\pi \in \mathcal{A}(\d, \pi_R, \delta) \\ \nu \in i\mathfrak{h}_M^\star}} \tau_{2}\left( \frac{\d}{\c(\pi^{R})} \right) \epsilon_\rho^{\delta, Q}(\tau, \nu) \dd\nu \right) \\
& \qquad + O\left(\sum_{\substack{\chi \in X^\mathrm{ur}(G) \\ \chi_R \simeq \pi_R}} m(\chi^{R}, \d) \widehat{\phi}(\chi)\right)
\end{align*}

\noindent that achieves the proof. \qed

\medskip

\noindent \textit{Remark.} Recall from Lemma \ref{lem:arch-tf-spectral} that the error term $\epsilon_\rho^{\delta, Q}$ is better when $\nu$ if far from the boundary of $\Omega$, so that it should be considered as a smoothed version of a bump function concentrating around the boundary, so that the integral \eqref{blouh} is a smoothed version of the counting number $B(\d, \pi_R, \delta, \partial_\rho Q; \phi)$, justifying the notation. This is jusified in \cite{brumley_counting_2018}.

\subsection{Towards the geometrical side}
\label{subsec:partition}

The equidistribution property has been recast as a convergence of spectral measures in Theorem \ref{thm:equid-prec}. The Selberg trace formula restates it as a geometrical quantity. In order to interpret  $B\left(\d,\pi_R, \delta, Q; \PPP\right)$ as a spectral side, it is necessary to add the contribution of the complementary part of the split archimedean spectrum. Summing the expressions above over all the spectral data and adding the complementary part of the spectrum,
\begin{equation*}
\begin{split}
\qquad \nu_Q(\widehat{\phi}\,) = 
& \frac{1}{Q^2}  \sum_{\substack{\substack{\sigma_R \in \widehat{G}_R \\ c(\sigma_R) \leqslant Q }}} \sum_{\substack{N\q \leqslant Q/c(\sigma_R)  \\ \q \wedge R = 1}}  \sum_{\substack{\underline{\delta} \in \mathcal{D} \\ \underline{\delta} = (M, \delta)}} \sum_{\d^S \div \q^S} \lambda_2 \left(\frac{\q^S}{\d^S}\right) J_\geom\left(\Phi_{\d, \pi_R, \delta, Q, \rho; \phi}\right) \\ 
& - \frac{1}{Q^2} \sum_{\substack{\substack{\sigma_R \in \widehat{G}_R \\ c(\sigma_R) \leqslant Q }}} \sum_{\substack{N\q \leqslant Q/c(\sigma_R)  \\ \q \wedge R = 1}}  \sum_{\substack{\underline{\delta} \in \mathcal{D} \\ \underline{\delta} = (M, \delta)}} \sum_{\d^S \div \q^S} \lambda_2 \left(\frac{\q^S}{\d^S}\right) J_\mathrm{comp}(\Phi_{\d, \pi_R, \delta, Q, \rho; \phi}) \\
& +  O \left(\frac{1}{Q^2}\sum_{\substack{\substack{\sigma_R \in \widehat{G}_R \\ c(\sigma_R) \leqslant Q }}} \sum_{\substack{N\q \leqslant Q/c(\sigma_R)  \\ \q \wedge R = 1}}  \sum_{\substack{\underline{\delta} \in \mathcal{D} \\ \underline{\delta} = (M, \delta)}} \sum_{\d^S \div \q^S} \lambda_2 \left(\frac{\q^S}{\d^S}\right) \Xi(\sigma_R, \phi) \right) \\
& +  O\left( \frac{1}{Q^2} \sum_{\substack{\substack{\sigma_R \in \widehat{G}_R \\ c(\sigma_R) \leqslant Q }}} \sum_{\substack{N\q \leqslant Q/c(\sigma_R)  \\ \q \wedge R = 1}}  \sum_{\substack{\underline{\delta} \in \mathcal{D} \\ \underline{\delta} = (M, \delta)}} \sum_{\d^S \div \q^S} \lambda_2 \left(\frac{\q^S}{\d^S}\right) \partial_\rho B(\d, \pi_R, \delta, Q) \right)
\label{decdec}
\end{split}
\end{equation*}

The main contribution  is carried by the first term, the remaining ones being showed below to contribute as negligible terms. Decompose the geometrical side $J_\geom(\Phi)$ as sum of two terms, the first one corresponding to the identity contribution\index{identity contribution, $\nu_{1, Q}$}, and the other being the elliptic remainder, in other words
\begin{equation}
J_{\geom}(\Phi) = \vol\left(G\left(F\right) \backslash G\left(\A\right)\right) \Phi(1)  + J_{\ell}(\Phi),
\end{equation}

\noindent where the elliptic part is expressed in term of orbital integrals
\begin{align*}
J_{\ell}\left(\Phi\right) & = \sum_{\{\gamma\} \neq \{1\}} \vol\left(G_\gamma(F) \backslash G_\gamma(\A)\right) \int_{G_\gamma(\A) \backslash G(\A)} \Phi\left(x^{-1}\gamma x\right) \dd x.
\end{align*}

\noindent The universal family counting measure now decomposes, via the splitting above, as
\begin{equation}
\nu_Q  = \vol\left(G\left(F\right) \backslash G\left(\A\right)\right) \nu_{1,Q} + \nu_{\ell, Q} - \nu_{\mathrm{comp}, Q} + O( \nu_{\Xi, Q}) + O(\nu_{\partial, Q}),
\label{splitting}
\end{equation}

\noindent where the following measures have been introduced.
\begin{align*}
\qquad \nu_{1,Q}\left(\widehat{\phi}\right) = &  \frac{1}{Q^2} \sum_{\substack{\substack{\sigma_R \in \widehat{G}_R \\ c(\sigma_R) \leqslant Q }}} \sum_{\substack{N\q \leqslant Q/c(\sigma_R)  \\ \q \wedge R = 1}}  \sum_{\substack{\underline{\delta} \in \mathcal{D} \\ \underline{\delta} = (M, \delta)}} \sum_{\d^S \div \q^S} \lambda_2 \left(\frac{\q^S}{\d^S}\right) {\Phi}_{\d, \pi_R, \delta, Q, \rho; \phi}(1) \\ 
\nu_{\ell , Q}\left(\widehat{\phi}\right)= & \frac{1}{Q^2} \sum_{\substack{\substack{\sigma_R \in \widehat{G}_R \\ c(\sigma_R) \leqslant Q }}} \sum_{\substack{N\q \leqslant Q/c(\sigma_R)  \\ \q \wedge R = 1}}  \sum_{\substack{\underline{\delta} \in \mathcal{D} \\ \underline{\delta} = (M, \delta)}} \sum_{\d^S \div \q^S} \lambda_2 \left(\frac{\q^S}{\d^S}\right) J_\mathrm{ell}({\Phi}_{\d, \pi_R, \delta, Q, \rho; \phi})\\
\nu_{\mathrm{comp},Q}\left(\widehat{\phi}\right) = &  \frac{1}{Q^2} \sum_{\substack{\substack{\sigma_R \in \widehat{G}_R \\ c(\sigma_R) \leqslant Q }}} \sum_{\substack{N\q \leqslant Q/c(\sigma_R)  \\ \q \wedge R = 1}}  \sum_{\substack{\underline{\delta} \in \mathcal{D} \\ \underline{\delta} = (M, \delta)}} \sum_{\d^S \div \q^S} \lambda_2 \left(\frac{\q^S}{\d^S}\right) J_\mathrm{comp}({\Phi}_{\d, \pi_R, \delta, Q, \rho; \phi}) \\
\nu_{\Xi,Q}\left(\widehat{\phi}\right)= &  \frac{1}{Q^2} \sum_{\substack{\substack{\sigma_R \in \widehat{G}_R \\ c(\sigma_R) \leqslant Q }}} \sum_{\substack{N\q \leqslant Q/c(\sigma_R)  \\ \q \wedge R = 1}}  \sum_{\substack{\underline{\delta} \in \mathcal{D} \\ \underline{\delta} = (M, \delta)}} \sum_{\d^S \div \q^S} \lambda_2 \left(\frac{\q^S}{\d^S}\right) \Xi(\sigma_R, \phi) \\ 
 \nu_{\partial,Q}\left(\widehat{\phi}\right)= &   \frac{1}{Q^2}\sum_{\substack{\substack{\sigma_R \in \widehat{G}_R \\ c(\sigma_R) \leqslant Q }}} \sum_{\substack{N\q \leqslant Q/c(\sigma_R)  \\ \q \wedge R = 1}}  \sum_{\substack{\underline{\delta} \in \mathcal{D} \\ \underline{\delta} = (M, \delta)}} \sum_{\d^S \div \q^S} \lambda_2 \left(\frac{\q^S}{\d^S}\right) \partial_\rho(\d, \sigma_R, \delta, Q)
\end{align*}

\section{Identity contribution}
\label{sec:identity}

For a given $\phi \in \mathcal{H}(G_S)$, the main term of $\nu_Q(\widehat{\phi}\,)$ is given by the contribution $\nu_{1,Q}(\widehat{\phi}\,)$ of the identity, and the other terms will be shown to be negligible. This section is dedicated to the computation of this identity contribution.

\begin{proposition}\label{prop:id} The contribution of the identity is, for $\phi \in \mathcal{H}(G_S)$,
\begin{equation*}
\vol\left(G\left(F\right) \backslash G\left(\A\right)\right) \nu_{1,Q}\left(\widehat{\phi}\right)= \nu\left(\widehat{\phi}\right)
+
\left\{
\begin{array}{cl}
O(Q^{-1} \log Q) & \text{ if $B$ is totally definite and $F=\Q$}\\
O(Q^{-\delta_F}) & \text{ if $B$ is totally definite and $F\neq \Q$}\\
O(\log^{-1} Q) & \text{ if there is a split infinite place} 
\end{array}
\right.
\end{equation*}

\noindent In particular, $\vol\left(G\left(F\right) \backslash G\left(\A\right)\right)  \nu_{1,Q}$ equidistributes with respect to $\nu$.
\end{proposition}


\subsection{Evaluating the test function at $1$}
\label{subsec:ft1}
Before summing over the spectral data, it is necessary to look at the inner part of $\nu_{1,Q}\left(\widehat{\phi}\right)$. Fix $\d$ an ideal of $\mathcal{O}^R$, $\pi_R$ a unitary irreducible representation of $G_R$ and $\delta$ a discrete archimedean parameter; and let for this section $\Phi = \Phi_{\d,\pi_R, \delta, Q, \rho;  \PPP}$ and $\Omega$ denote $\Omega(Q/N\d^S c(\sigma_R))$ for convenience. The very definition \eqref{def:ft} of $\Phi$ gives
\begin{equation}	
\Phi(1) = \varepsilon_{K_0(\d^S)}(1)  \phi_{S, f}^R(1) \xi_{\sigma_{R}} (1) \widehat{\phi}_R(\pi_R) f_\rho^{\delta, Q, \phi}(1).
\end{equation}

\subsubsection{Finite split places out of ${S}$}

\noindent For the prime-to-$S$ split finite part, by definition
\begin{equation}
\varepsilon_{\overline{K}_{0}\left(\d^S\right)}\left(1\right)  = \vol\left(\overline{K}_0(\d^S)\right)^{-1}.
\end{equation}

\noindent The volume of a cofinite subgroup depends on its index, and the indices of classical congruence subgroups are well-known \cite{diamond_first_2005}. Introduce $K^{R,S} = \prod_{v \notin R\cup S} K_v$. Since $Z^{R,S}$ is fully contained in $K^{R,S}_0(\d^S)$ for all ideal $\d^S$, 
\begin{equation}
\label{sgc-volume}
\left[\overline{K}^{R,S}:\overline{K}_0\left(\d^S\right)\right] = \left[K^{R,S}:K_0\left(\d^S\right)\right],
\end{equation}

\noindent by the isomorphism theorems. So thanks to the normalizations chosen for the measures,
\begin{equation}
\varepsilon_{\overline{K}_{0}\left(\d^S\right)}\left(1\right)  =  \left[K^{R,S}:K_0\left(\d^S\right)\right] = \left(\id \star \mu^2\right) (\d^S) =: \varphi_2(\d^S).
\label{volume}
\end{equation}

\subsubsection{Finite split places in $S$}

\noindent For the $S$-split finite part, the Plancherel inversion formula \eqref{plancherel} gives
\begin{equation}
\phi_{S, f}^R(1) = \int_{\widehat{G}_{S, f}^R} \widehat{\phi}_{S, f}^R (\pi_{S, f}^R) \dd \pi_{S, f}^R.
\end{equation}

\subsubsection{Ramified places}

\noindent For the ramified matrix coefficient \eqref{matrix-coeff:selecting-function}, by the Plancherel formula \eqref{plancherel} and the normalization chosen for $\xi_{\pi_R}$,
\begin{equation}
\xi_{\pi_R}\left(1\right) = \int_{\widehat{G}_R} \mathbf{1}_{\sigma \simeq \pi_R} \dd\mu_R^\Pl( \sigma) =  \mu^{\Pl}_R(\pi_R).
\end{equation}

\subsubsection{Split archimedean places}

The test function $f_\rho^{\delta, Q, \phi}$ at archimedean split places is not so immediate to evaluate at $1$, since it is not an explicit function but provided by an existence theorem. The Plancherel formula allows to express it in terms of its Fourier transform $h_\rho^{\delta, Q, \phi}$ for which Lemma \ref{lem:arch-tf-spectral} provides information. There is an error term function due to the smoothing procedure willingly kept undisclosed and for which bounds are provided later.

\begin{lemma}
\label{lem:arch-tf-geometric}
For every $\delta \in \mathcal{D}$, $\rho > 0$, $\phi \in \mathcal{H}(G)$ and bounded set of tempered parameters $\Omega$,
\begin{equation}
f_\rho^{\delta, Q, \phi}(1) = \int_{\Omega} \widehat{\phi}(\pi_{{\delta, \nu}}) \dd\nu + \partial_\rho B(\delta, Q),
\end{equation}

\noindent where
\begin{equation}
\partial_\rho B(\delta, Q) = \int_{i\h_M^\star} \varepsilon_\rho^{\delta, Q}(\pi_{\delta, \nu}) \dd \nu.
\end{equation}
\end{lemma}

\proof This is the proof of \cite[Lemma 11.2]{brumley_counting_2018} into which the spectral localizing function around $(\delta, Q)$ they build, namely $h_\rho^{\delta, Q}$ in their notation for reference, is replaced by $h_\rho^{\delta, Q, \phi}$. The Plancherel formula gives
\begin{equation}
f_\rho^{\delta, Q, \phi}(1) = \int_{i\h_M^\star} h_\rho^{\delta, Q}(\pi_{{\delta, \nu}}) \widehat{\phi}(\pi_{{\delta, \nu}}) \dd \nu.
\end{equation}

Integrating the approximation of Lemma \ref{lem:arch-tf-spectral} yields
\begin{equation}
\label{archimedean-selecting-approx}
\int_{i\h_M^\star} h_\rho^{\delta, Q}(\pi_{\delta, Q}) \widehat{\phi}(\pi_{{\delta, \nu}}) \dd \nu =  \int_{\Omega} \widehat{\phi}(\pi_{{\delta, \nu}}) \dd \nu + \int_{i\h_M^\star} \varepsilon_\rho^{\delta, Q}(\pi_{\delta, \nu}) \dd \nu ,
\end{equation}

\noindent and this achieves the proof. \qed

\medskip

Finally it follows a more explicit form of the value at the identity, namely
\begin{equation}
\label{value1}
\begin{split}
\Phi\left(1\right) & = \varphi_2 (\d^S)  \mu^{\Pl}_R(\pi_R) \widehat{\phi}(\pi_R) \int_{\Omega} \widehat{\phi}(\pi_{\delta, \nu})\dd\nu\int_{\substack{\pi_{S,f}^R \in \widehat{G}_{S,f}^R}} \widehat{\phi}(\pi_{S,f}^R)\dd\pi_{S,f}^R \\
& \qquad + \varphi_2 (\d^S)  \mu^{\Pl}_R(\pi_R) \widehat{\phi}(\pi_R) \partial_\rho B(\delta, Q) \int_{\substack{\pi_{S,f}^R \in \widehat{G}_{S,f}^R}} \widehat{\phi}(\pi_{S,f}^R)\dd\pi_{S,f}^R 
\end{split}
\end{equation}

The tools are now in place to, after summation of \eqref{value1}, provide a decomposition of the counting measure of the whole universal family.

\subsection{Splitting the identity contribution}
\label{subsec:estimations1}

Recall that $\phi$ fixes the $S$-part of the conductor, so that every $S$-part of ideal appearing from now on is fixed, namely the only one non killed by the action of $\phi$. However, the choice made is to keep formulations in terms if ideals of the whole integer ring $\mathcal{O}$, as more convenient and helping to think about the counting law situation where $S$ is empty. The following decomposition holds for the identity part of the counting measure.
\begin{proposition}
\label{prop:nuQ1}
For every $Q \geqslant 1$, 
\begin{equation}
\nu_{1,Q} = \nu_{1,Q}^{(p)}+\nu_{1,Q}^{(e1)}+\nu_{1,Q}^{(e2)},
\end{equation}

\noindent where $\nu_{1,Q}^{(p)}$ is the main identity term, namely
\begin{equation}
\label{nu-p}
\begin{split}
\nu_{1,Q}^{(p)}(\widehat{\phi}\,) & =  \frac{1}{2}\frac{\zeta^{S,R \star}(1)\zeta^{S,R}(2)}{\zeta^{S,R}(4)} \int_{\substack{\pi_{S,f}^R \in \widehat{G}_S^R}} \frac{\widehat{\phi}(\pi_{S,f}^R)}{c(\pi^R_S)^2} \sum_{\substack{N\m^S \leqslant Q/c(\pi^R_S) \\ \m^S \wedge R = 1}} \frac{\lambda_2(\m^S)}{(N\m^S)^2}    \\
&\sum_{\substack{\pi_R \in \widehat{G}_R \\ c\left(\pi_R\right) \leqslant Q/N\m^S c(\pi^R_S)}}  \frac{\widehat{\phi} (\pi_R)}{c(\pi_R)^2}\mu^{\Pl}_R(\pi_R)  \sum_{\substack{\underline{\delta} \in \mathcal{D} \\ \underline{\delta} = (M, \delta)}} \int_{\Omega_{\delta}(Q/N\m^S c(\pi_R) c(\pi^R_S))} \frac{\widehat{\phi}(\pi_{\delta, \nu})}{c(\pi_{\delta, \nu})^{2}}\dd\nu \dd\pi_{S,f}^R 
\end{split}
\end{equation}

\noindent and $\nu_{1,Q}^{(e1)}(\widehat{\phi}\,)$ is the error term due to the smoothing, namely
\begin{equation}
\label{nu-e1}
\begin{split}
\nu_{1,Q}^{(e1)} \left(\widehat{\phi}\right)  & \ll  
\frac{1}{Q^2}  \int_{\substack{\pi_{S,f}^R \in \widehat{G}_S^R }} \widehat{\phi}(\pi_{S,f}^R)  \\ 
& \qquad  \sum_{\substack{N\q \leqslant Q/c(\pi^R_S)\\\q \wedge R = 1}} \sum_{\d^S \div \q^S}   \lambda_2\left(\frac{\q^S}{\d^S}\right)\varphi_2(\d^S)  \sum_{\substack{\pi_R \in \widehat{G}_R \\ c\left(\pi_R\right) \leqslant Q/N\q^S c(\pi^R_S)}} \mu^{\Pl}_R(\pi_R) |\widehat{\phi}| (\pi_R) \\ 
& \qquad \sum_{\substack{\underline{\delta} \in \mathcal{D} \\ \underline{\delta} = (M, \delta)}} \partial_\rho B(\delta, Q(Q/N\m^S c(\pi_R) c(\pi_S^R)))  \dd\pi_{S,f}^R
\end{split}
\end{equation}

\noindent and $\nu_{1,Q}^{(e2)}(\widehat{\phi}\,)$ is an extra error term, that is

\begin{equation}
\label{nu-e2}
\begin{split}
\nu_{1,Q}^{(e2)}(\widehat{\phi}\,) & \ll  Q^{-\delta_F+\varepsilon_F} \int_{\substack{\pi_{S,f}^R \in \widehat{G}_S^R}} \frac{\widehat{\phi}(\pi_{S,f}^R)}{c(\pi^R_S)^{2-\delta_F+\varepsilon_F}}   \sum_{\substack{N\m^S \leqslant Q/c(\pi^R_S) \\ m^S \wedge R = 1}} \frac{\lambda_2(\m^S)}{(N\m^S)^{2-\delta_F+\varepsilon_F}}  \\
& \qquad \sum_{\substack{\pi_R \in \widehat{G}_R \\ c\left(\pi_R\right) \leqslant Q/N\m^S c(\pi^R_S)}} \frac{\widehat{\phi} (\pi_R)}{c(\pi_R)^{2-\delta_F+\varepsilon_F}}\mu^{\Pl}_R(\pi_R)  \\
& \qquad  \sum_{\substack{\underline{\delta} \in \mathcal{D} \\ \underline{\delta} = (M, \delta)}} \int_{\Omega_{\delta}(Q/N\m^S c(\pi_R)c(\pi^R_S))} \frac{\widehat{\phi}(\pi_{\delta, \nu})}{c(\pi_{\delta, \nu})^{2-\delta_F+\varepsilon_F}}\dd\nu \dd\pi_{S,f}^R
\end{split}
\end{equation}
\end{proposition}

\proof The counting measure has been decomposed in measures on harmonic subfamilies \eqref{step3} of fixed spectral parameters. These measures have been given a geometric interpretation by the mean of the trace formula in Lemma \ref{lem:selecting}, whose identity contribution \eqref{value1} is given above. After summation of the identity contributions over the spectral data constituting the truncated universal family,
\begin{equation*}
\label{identity-splitting}
\begin{split}
\nu_{1,Q} \left(\widehat{\phi}\right)  & =  
\frac{1}{Q^2} \int_{\substack{\pi_{S,f}^R \in \widehat{G}_{S,f}^R}}  \widehat{\phi}(\pi_{S,f}^R) \sum_{\substack{N\q \leqslant Q/c(\pi_{S,f}^R) \\\q \wedge R = 1}} \sum_{\d^S \div \q^S}   \lambda_2\left(\frac{\q^S}{\d^S}\right)\varphi_2(\d^S)    \\
& \sum_{\substack{\pi_R \in \widehat{G}_R \\ c\left(\pi_R\right) \leqslant Q/N\d c(\pi_{S,f}^R)}} \mu^{\Pl}_R(\pi_R) \widehat{\phi} (\pi_R) \sum_{\substack{\underline{\delta} \in \mathcal{D} \\ \underline{\delta} = (M, \delta)}}  \int_{\Omega_\delta(Q/N\d c(\pi_R)\pi_{S,f}^R )} \widehat{\phi}(\pi_{\delta, \nu})\dd\nu \dd\pi_{S,f}^R .
\end{split}
\end{equation*}

Sums of arithmetic functions on ideals of number fields can be explicitly evaluated. This motivates a permutation of sums and integrals in order to estimate the sum over the volumes $\varphi_2(\d^S)$  first, so that
\begin{align*}
\nu_{1,Q} \left(\widehat{\phi}\right)  & = \frac{1}{Q^2} \int_{\substack{\pi_{S,f}^R \in \widehat{G}_S^R}} \widehat{\phi}(\pi_{S,f}^R) \sum_{\substack{N\m^S \leqslant Q/c(\pi^R_{S, f}) \\ \m^S \wedge R = 1}} \lambda_2(\m^S) \sum_{\substack{\pi_R \in \widehat{G}_R \\ c\left(\pi_R\right) \leqslant Q/N\m}} \mu^{\Pl}_R(\pi_R) \widehat{\phi} (\pi_R) \\ 
& \qquad \sum_{\substack{\underline{\delta} \in \mathcal{D} \\ \underline{\delta} = (M, \delta)}} \int_{\Omega_{\delta}(Q/N\m c(\pi_R))} \widehat{\phi}(\nu) 
   \sum_{\substack{N\d^S \leqslant Q/ N\m c(\pi_R)c(\pi_{\delta, \nu})\\ \d \wedge R = 1}}\varphi_2(\d^S) \dd\nu \dd\pi_{S,f}^R  .
\end{align*}

The following lemma estimates the innermost sum.
\begin{lemma}
\label{lemma}
Let $\zeta^{S, R}$ be the prime-to-$R$-and-$S$ part of the zeta function associated to $F$, and $\zeta^{S, R\star}(1)$ its residue at $1$. For any $X>0$, 
\begin{equation}
\sum_{\substack{N\d^S \leqslant X \\ \d \wedge R = 1}}  \varphi_2 (\d^S)  = \frac{1}{2}\frac{\zeta^{S, R \star}(1)\zeta^{S,R}(2)}{\zeta^{S,R}(4)} X^2   + 
\left\{
\begin{array}{cl}
O(X \log X) & \text{ if } F = \Q, \\
O(X^{2-\delta_F}) & \text{ otherwise}
\end{array}
\right. 
\end{equation}
\end{lemma}

\rk It is possible to note a posteriori that the remainder term shown here is sharp, and it gives rise to the most significant remainder appearing in Theorem \ref{thm-count} and Theorem \ref{thm:equid-prec}, provided $F \neq \Q$, except the one coming from the smoothing part detailed in Lemma \ref{lem:arch-tf-geometric} that is absent from the totally definite setting. Hence, provided the smoothing problem can be solved and thus a sharp count realized without excessive loss, the error would have power savings and will be similar to the totally definite case. 

\medskip

\proof Remind that all the ideals superscrited $S$ are prime to $S$. Standard estimates of the sum of ideals given by \cite{landau_einfuhrung_1918} lead to
\begin{align*}
\sum_{\substack{N\d^S \leqslant X \\ \d \wedge R = 1}}   \varphi_2 (\d^S) &   = \sum_{\substack{N\l^S\leqslant X\\ \l^S \wedge R = 1}}  \mu^2(\l^S) \sum_{\substack{N\m^S \leqslant X/N\l \\ \m^S \wedge R = 1}} N\m^S  \\
& = \sum_{\substack{N\l^S\leqslant X \\ \l^S \wedge R = 1}}  \mu^2(\l^S) \left[ \frac{\zeta^{S, R \star}(1)}{2} \frac{X^2}{(N\l^S)^2} + O\left(\left(\frac{X}{N\l^S}\right)^{2-\delta_F}\right) \right] \\
& = \frac{1}{2}\zeta^{S, R \star}(1) X^2  \sum_{\substack{N\l^S\leqslant X\\ \d \wedge R = 1}}  \frac{\mu^2(\l^S)}{(N\l^S)^2}  + O\left(X^{2-\delta_F} \sum_{\substack{N\l^S\leqslant X\\ \l^S \wedge R = 1}}  \frac{\mu^2(\l^S)}{(N\l^S)^{2-\delta_F}}\right) \\
& = \frac{1}{2}\frac{\zeta^{S,R \star}(1)\zeta^{S,R}(2)}{\zeta^{S,R}(4)} X^2   + 
\left\{
\begin{array}{cl}
O(X \log X) & \text{ if } F = \Q;  \\
O(X^{2-\delta_F}) & { otherwise};
\end{array}
\right. 
\label{ev1}
\end{align*}

\noindent where the knowledge of the Dirichlet series associated to $\mu^2$ yielded
\begin{equation}
\sum_{N(\m) \leqslant X} \frac{\mu^2(\m)}{N\m} \sim \frac{\zeta^\star(1)}{\zeta(2)}\log X = O(\log X),
\end{equation}

\noindent in the case $F=\Q$, giving the worst remainder term. Otherwise, the sum is convergent. \qed

\medskip

This lemma induces a splitting of $\nu_{1,Q}$ as $\nu_{1,Q}^{(p)} + \nu_{1,Q}^{(e2)}$ according to the principal and error parts in the lemma above and, adding the error term coming from the smoothing evaluation at the identity \eqref{value1}, achieves to prove the claim. \qed

\subsection{Estimating the main part $\nu_{1,Q}^{(p)}$}
\begin{proposition}
\label{prop:main-part-id-2}
For every $Q \geqslant 1$, the main part admits the asymptotic development
\begin{equation}
\label{prop:main-part-id}
\vol\left(G\left(F\right) \backslash G\left(\A\right)\right) \nu^{(p)}_{1,Q}\left(\widehat{\phi}\,\right)= \nu\left(\widehat{\phi}\,\right) + O(Q^{-2}).
\end{equation}
\end{proposition}

\proof Recall the term $\nu_{1,Q}^{(p)}$ of Proposition \ref{prop:nuQ1}, namely
\begin{align*}
\nu_{1,Q}^{(p)}(\phi) & =  \frac{1}{2}\frac{\zeta^{S,R \star}(1)\zeta^{S,R}(2)}{\zeta^{S,R}(4)} \int_{\substack{\pi_{S,f}^R \in \widehat{G}_S^R}} \frac{\widehat{\phi}(\pi_{S,f}^R)}{c(\pi^R_S)^2}    \\
& \qquad \sum_{\substack{N\m^S \leqslant Q/c(\pi^R_S) \\ \m^S \wedge R = 1}} \frac{\lambda_2(\m^S)}{(N\m^S)^2}  \sum_{\substack{\pi_R \in \widehat{G}_R \\ c\left(\pi_R\right) \leqslant Q/N\m^S c(\pi^R_S)}} \frac{\widehat{\phi} (\pi_R)}{c(\pi_R)^2}\mu^{\Pl}_R(\pi_R) \\
& \qquad  \sum_{\substack{\underline{\delta} \in \mathcal{D} \\ \underline{\delta} = (M, \delta)}} \int_{\Omega_{\delta}(Q/N\m^S c(\pi_R)c(\pi_S^R))} \frac{\widehat{\phi}(\pi_{\delta, \nu})}{c(\pi_{\delta, \nu})^{2}}\dd\nu  \dd\pi_{S,f}^R 
\end{align*}

\noindent The following lemmata state the convergence of the integral over archimedean parameters and of the sum over ramified parts.
\begin{lemma}
\label{summeasures}
For every $\Re(s) >1$, the following sum converges as $Q \to \infty$.
\begin{equation}
\sum_{\substack{\pi_R \in \widehat{G}_R \\ c(\pi_R) \leqslant Q}} \frac{\mu^{\Pl}_R(\pi_R)}{c(\pi_R)^s}.
\end{equation}
\end{lemma}

\proof The Jacquet-Langlands correspondence states a bijection between $\widehat{G}_R$ and the discrete part of the spectrum of $\widehat{\PGL}(2,F_R)$, which preserves both formal degrees, which are the Plancherel measures $\mu_R^\Pl(\pi_R)$, and conductors by definition. Hence,
\begin{equation}
\sum_{\substack{\sigma_R \in \widehat{G}_R \\ c(\sigma_R) \leqslant Q}} \frac{\mu^{\Pl}_R(\sigma_R)}{c(\sigma_R)^s}  \leqslant \sum_{\pi_R \in \widehat{\PGL}(2,F_R)^\mathrm{disc}} \frac{\mu^{\Pl}_R(\pi_R)}{c(\pi_R)^s},
\end{equation}

\noindent and that last sum is finite for $\Re(s)>1$ by the case of $\PGL(2)$ by the computations of \cite{brumley_counting_2018} or by Section \ref{subsec:constant} below. Hence, it follows the sought convergence for the ramified parts, ending the proof of the lemma. $\quad \square$

\begin{lemma}
\label{summeasures2}
\label{lem:plancherel-cv}
For every $\Re(s) >1$, the following integral converges absolutely as $X \to \infty$.
\begin{equation}
\int_{\Omega_{\delta}(X)} \frac{\widehat{\phi}(\pi_{\delta, \nu})}{c(\pi_{\delta, \nu})^s} \dd\nu.
\end{equation}
\end{lemma}

This is Lemma 6.12 of \cite{brumley_counting_2018} concerning the $\GL(2)$ case.  Let us denote $\displaystyle \int_{\widehat{G}_R} \frac{\widehat{\phi}(\pi_R)}{c(\pi_R)^s}\dd\pi_R$ and $\displaystyle \int_{\widehat{G}^R_\infty} \frac{\widehat{\phi}(\pi_R)}{c(\pi^R_\infty)^s}\dd\pi^R_\infty$ the limits in the above lemmata. The prime-to-$S$-and-$R$ part of the Dirichlet series associated to $\lambda_2$ converges at $2$ to $\zeta^{S,R}_F(2)^{-2}$ and makes the expression of $\nu_{1, Q}^{(p)}$ converges to
\begin{equation}
\label{mt}
\frac{1}{2}\frac{\zeta^{S,R \star}(1)}{\zeta^{S,R}(2)\zeta^{S,R}(4)}\int_{\widehat{G}_{R}} \frac{\widehat{\phi}(\pi_R)}{c(\pi_{R})^2}\dd\pi_{R}\int_{\widehat{G}^R_S} \frac{\widehat{\phi}(\pi_S^R)}{c(\pi_{S}^R)^2}\dd\pi_{S}^R\int_{\widehat{G}_{\infty}^R} \frac{\widehat{\phi}(\pi_\infty^R)}{c(\pi_{\infty}^R)^2}\dd\pi_{\infty}^R.
\end{equation}

\subsection{Rewriting the constant}
\label{subsec:constant}

Previous computations unveiled the constant
\begin{equation}
\frac{\zeta^{S, R \star}(1)}{\zeta^{S,R}(2)\zeta^{S,R}(4)}.
\end{equation}

\noindent It is possible to give to this constant a more geometrical flavour by reformulating the special values of the zeta functions appearing in terms of volumes. This is the content of the following lemma.
\begin{proposition}
For every finite set of places $S$, 
\begin{equation}
\frac{\zeta^{S, R \star}(1)}{\zeta^{S,R}(2)\zeta^{S,R}(4)} = \int_{\widehat{G}^{S,R}}^\star \frac{\dd\pi^{S,R}}{c(\pi^{S,R})^2} = \zeta^{S,R \star}(1) \prod_{\p \notin S \cup R} \zeta_\p(1)^{-1} .
\end{equation}
\end{proposition}

\proof The knowledge of the volumes of congruence subgroups \eqref{volume} gives
\begin{equation}
\varepsilon_{\overline{K}_{0, \p}(\p^r)} (1) = \vol\left( \overline{K}_{0, \p}(\p^r)\right)^{-1} = (\id \star \mu^2)(\p^r).
\end{equation}

On an other hand, this volume can be computed by the Plancherel formula. Introduce the volume of slices of the spectrum of fixed conductor
\[
M_\p(\p^r) = \int_{\substack{\sigma_\p \in \widehat{G}_\p \\ \c(\sigma_\p) = \p^r}} \dd\sigma_\p, \qquad r \geqslant 1. .
\] 

\noindent The Plancherel inversion formula then yields
\begin{align*}
\varepsilon_{\overline{K}_{0, \p}(\p^r)}(1) & = \int_{\widehat{G}_\p} \widehat{\varepsilon}_{\overline{K}_{0, \p}(\p^r)}(\pi_\p)\dd\pi_\p  = \int_{\widehat{G}_\p} \tau_2\left(\frac{\p^r}{\c(\pi_\p)}\right) \dd \pi_\p \\
& = \sum_{\d \div \p^r} M_\p(\d) \tau_2\left(\frac{\p^r}{\d}\right) = (M_\p \star \tau_2)(\p^r)
\end{align*}

\noindent Hence, by inversion, $M_\p = \id \star \mu^2 \star \lambda_2$. In particular, the local Dirichlet series associated to $M_\p$ is given by
\begin{equation}
D_\p\left(s\right) = \sum_{\substack{\m = \p^r \\ r \geqslant 0}} \frac{M_\p(\m)}{N\m^s} = \frac{\zeta_\p(s-1)}{\zeta_\p(s)\zeta_\p(2s)}, \qquad \Re(s) > 1.
\end{equation}

\noindent Evaluating it at $s=2$, a new expression for the local special values appearing in the constant is
\begin{equation}
\label{special-values-local-measures}
\int_{\widehat{G}_\p} \frac{\dd\pi_\p}{c(\pi_\p)^2}  = \frac{\zeta_\p(1)}{\zeta_\p(2)\zeta_\p(4)},
\end{equation}

\noindent proving the finiteness of the local integrals defining the equidistribution measure \eqref{mu} at the finite places, as claimed in the introduction. However, the infinite product over $\p \notin R$ of these quantities unfortunately diverges, for $1$ is a pole of $\zeta^{S,R}$. This motivates a slight modification in order to compensate it by the residue at 1. Introduce the regularized integral
\begin{equation*}
\int_{\widehat{G}^{S,R}}^\star \frac{\dd\pi^{S,R}}{c(\pi^{S,R})^2} = \zeta^{S,R \star}(1) \prod_{\p \notin S \cup R} \zeta_\p(1)^{-1}  \int_{\widehat{G}_\p} \frac{\dd\pi_\p}{c(\pi_\p)^2} = \zeta^{S,R \star}(1)\prod_{\p \notin S \cup R} \frac{1}{\zeta_\p(2)\zeta_\p(4)}, 
\end{equation*}

\noindent ending the proof. \qed

\medskip

%

\noindent The global integral\index{regularized integral} is defined to be
\begin{equation*}
\int_{\widehat{\Pi}}^\star \frac{\dd\pi}{c(\pi)^2} = \int_{\widehat{G}^{S,R}}^\star \frac{\dd\pi^{S,R}}{c(\pi^{S,R})^2} \int_{\widehat{G}_{S \cup R}} \frac{\widehat{\phi}(\pi_{S \cup R})}{c(\pi_{S \cup R})^2} \dd\pi_{S \cup R} =  \zeta^{\star}(1) \prod_{v} \zeta_v(1)^{-1}  \int_{\widehat{G}_v} \frac{\dd\pi_v}{c(\pi_v)^2}.
\end{equation*}

\noindent It thus follows the  expression \eqref{mu} of the regularized integral, giving the desired statement and motivating the choice of both the measure $\mu$ and the constant $C$. Since the error terms in Lemmata \ref{summeasures} and \ref{summeasures2} are those of Dirichlet series at a point distant by 1 from their abscissa of convergence, the expression \eqref{mt} rewrites
\begin{equation}
\nu_{1,Q}^{(p)}(\widehat{\phi}\,) = \frac{1}{2}   \int_{\widehat{\Pi}}^\star \frac{\widehat{\PPP}(\pi)}{c(\pi)^2} \dd\pi  + 
O(Q^{-1}),
\end{equation}

\noindent reaching the term of the proof of Proposition \ref{prop:main-part-id-2}. \qed

\medskip

\begin{rk}
Notice that the Sauvageot theorem is a two-edged result: it opens the path to equidistribution and allows conclusions for characteristic functions which are not of the form $\widehat{\phi}$; however it also spoils the remainder term for general functions. This error term remains only for specific functions either admissible, \textit{i.e.} of the form $\widehat{\phi}$ for $\phi$ in the Hecke algebra of $G$, in particular for the counting problem.
\end{rk}

\subsection{Estimating the smoothing part $\nu_{1,Q}^{(e1)}$}
\label{subsec:smoothing-error}

The following lemma states the negligibility of the smoothing part compared to the main one.

\begin{lemma}
\label{lem:smoothing-part}
For every $Q \geqslant 1$, 
\begin{equation}
\nu_{1,Q}^{(e1)} \left(\, \widehat{\phi}\,\right)  \ll \log(Q)^{-1}.
\end{equation}
\end{lemma}

\proof The term \eqref{nu-e1} coming from the smoothed selection of archimedean parts is
\begin{equation}
\begin{split}
\nu_{1,Q}^{(e1)} \left(\widehat{\phi}\right)  & \ll  
\frac{1}{Q^2}  \int_{\substack{\pi_{S,f}^R \in \widehat{G}_S^R }} \widehat{\phi}(\pi_{S,f}^R)  \\ 
& \qquad  \sum_{\substack{N\q \leqslant Q/c(\pi^R_S)\\\q \wedge R = 1}} \sum_{\d^S \div \q^S}   \lambda_2\left(\frac{\q^S}{\d^S}\right)\varphi_2(\d^S)  \sum_{\substack{\pi_R \in \widehat{G}_R \\ c\left(\pi_R\right) \leqslant Q/N\q^S c(\pi^R_S)}} \mu^{\Pl}_R(\pi_R) \widehat{\phi} (\pi_R) \\ 
& \qquad \sum_{\substack{\underline{\delta} \in \mathcal{D} \\ \underline{\delta} = (M, \delta)}} \partial_\rho B(\delta, Q(Q/N\d^S c(\pi_R) c(\pi_S^R)))  \dd\pi_{S,f}^R, 
\end{split}
\end{equation}

The following lemma is a straighforward adaptation of the work of Brumley and Mili\'{c}evi\'{c} in which the contribution of the smoothing effect is negligible.

\begin{lemma}
\label{lem:BM-bounds}
The contribution of the smoothing error term satisfies, for a suitable choice of $\rho$ depending on $\q$,
\begin{align*}
\sum_{\substack{N\q \leqslant Q \\ \q \wedge R, S = 1}} & \sum_{\substack{\sigma_R \in \widehat{G}_R \\ c(\sigma_R) \leqslant Q/N\q}} \sum_{\delta \in\mathcal{D}} \sum_{\d | \q} \lambda_2\left(\frac{\q^S}{\d^S} \right) \phi_2(\d^S)\mu_R^\Pl(\sigma_R) \partial_\rho B(\delta, Q)  \ll \frac{Q^2}{\log(Q)}
\end{align*}
\end{lemma}

\proof This is in essence Proposition 12.1 in \cite{brumley_counting_2018}. Their Lemma 12.2 states
\begin{equation}
\sum_{\substack{N\q \leqslant Q \\ \q \wedge R, S = 1}} \sum_{\delta \in\mathcal{D}} N\q   \partial_\rho B(\delta, Q(Q/N\q)) \ll \frac{Q^2}{\log(Q)}.
\end{equation}

Two slight modifications have to be mentioned because of the existence of ramified places and of the use of a different filtration. The centerless setting leads to use a filtration of subgroups $K_{0}(\q)$ whose indices are of order $N\q$, instead of $N\q^2$ in their case, justifying the presence of $N\q$ in the equation above. The existence of ramified places is dealt with by plugging this estimates above in the whole sum, 
\begin{align*}
\sum_{\substack{N\q \leqslant Q \\ \q \wedge R, S = 1}} &  \sum_{\substack{\sigma_R \in \widehat{G}_R \\ c(\sigma_R) \leqslant Q/N\q}} \sum_{\delta \in\mathcal{D}} \sum_{\d | \q} \lambda_2\left(\frac{\q^S}{\d^S} \right) \phi_2(\q) \mu_R^\Pl(\pi_R) \partial_\rho B(\delta, Q(Q/N\q c(\pi_R)))   \\
& \ll \frac{Q^2}{\log Q}\sum_{\substack{\sigma_R \in \widehat{G}_R \\ c(\sigma_R) \leqslant Q/N\q}} \frac{\mu^\Pl_R(\pi_R)}{c(\pi_R)^{2-\varepsilon}}
\end{align*}

\noindent and this last sum converges by Lemma \ref{summeasures}, finishing the proof. \qed

\medskip

\rk Note that the worst error term in Theorems \ref{thm-count} and \ref{thm:equid-prec} comes from this part, corresponding to the smoothing of the selecting function for split archimedean places. As in the $\GL(2)$ case \cite{brumley_counting_2018}, this is where the desired power savings is lost, safe in the totally definite case when there is no involved smoothing.

\subsection{Estimating the error term $\nu_{1,Q}^{(e2)}$}

\begin{lemma}
\label{lem:extra-error-term}
For every $Q \geqslant 1$, 
\begin{equation}
\nu_{1,Q}^{(e2)}(\phi)\ll  Q^{-\delta_F+\varepsilon_F},
\end{equation}

\noindent for all $\varepsilon_F>0$ for the case $F=\Q$, and for $\varepsilon_F=0$ otherwise.
\end{lemma}

\proof Now turn back to treatment of the $\nu_{1,Q}^{(e2)}$ term coming from the remainder in Lemma \ref{lemma}. The bound that has to be refined is
\begin{align*}
\nu_{1,Q}^{(e2)}(\phi) & \ll  Q^{-\delta_F+\varepsilon_F} \int_{\substack{\pi_{S,f}^R \in \widehat{G}_S^R}} \frac{|\widehat{\phi}|(\pi_{S,f}^R)}{c(\pi^R_S)^{2-\delta_F+\varepsilon_F}}   \\
&   \qquad  \sum_{\substack{N\m^S \leqslant Q/c(\pi^R_S) \\ m^S \wedge R = 1}} \frac{\lambda_2(\m^S)}{(N\m^S)^{2-\delta_F+\varepsilon_F}}   \sum_{\substack{\pi_R \in \widehat{G}_R \\ c\left(\pi_R\right) \leqslant Q/N\m^S c(\pi^R_S)}} \frac{|\widehat{\phi}| (\pi_R)}{c(\pi_R)^{2-\delta_F+\varepsilon_F}}\mu^{\Pl}_R(\pi_R)   \\
& \qquad \sum_{\substack{\underline{\delta} \in \mathcal{D} \\ \underline{\delta} = (M, \delta)}} \int_{\Omega_{\delta}(Q/N\m^S c(\pi_R)c(\pi_S^R))} \frac{|\widehat{\phi}|(\pi_{\delta, \nu})}{c(\pi_{\delta, \nu})^{2-\delta_F+\varepsilon_F}}\dd\nu\dd\pi_{S,f}^R
\end{align*}

\noindent The inner sums and integrals converge by Lemmata \ref{summeasures} and \ref{summeasures2}, since $2-\delta_F+\varepsilon_F$ is always greater than $1$. It follows a remainder term in $Q^{-\delta_F+\varepsilon_F}$. \qed

\medskip

At last, the asymptotic development obtained in Proposition \ref{prop:main-part-id-2} and the bounds obtained in Lemmata \ref{lem:smoothing-part} and \ref{lem:extra-error-term} prove the equidistribution of the identity part of the counting measure with respect to $\nu$, as stated in Proposition \ref{prop:id}.

\section{Spectral error terms}
\label{sec:REMAINDER}

\subsection{Characters contribution}
\label{subsec:characters}

Recall that the global characters contribution is given by
\begin{equation}
\label{charac_contrib}
\nu_{\Xi, Q}(\widehat{\PPP}) = \frac{1}{Q^2}\sum_{\substack{N\q \leqslant Q \\ \q \wedge R = 1}} \sum_{\substack{\pi_R \in \widehat{G}_R \\ c(\pi_R) \leqslant Q/N\q}} \sum_{\d^S \div \q^S} \lambda_2\left(\frac{\q^S}{\d^S}\right) \Xi(\PPP, \pi_R).
\end{equation}

\begin{lemma}
\label{characrers}
\label{lem:characrers}
For every $\varepsilon > 0$,
\begin{equation}
\nu_{\Xi, Q}(\widehat{\PPP}) \ll Q^{-1+\varepsilon}.
\end{equation}
\end{lemma}

\proof Similarly to the intervention of the trace formula to make explicit the measure $\nu_Q$, the Poisson summation formula is the main tool to count characters. The counting measure for characters can be interpreted as a spectral side, such that every non-identity terms vanishes on the geometric side. Recall that for a character $\pi_R$, since the multiplicities are all equal to one,
\begin{equation}
\Xi(\pi_R, \PPP) = \sum_{\substack{\chi \in X^{\mathrm{ur}}(G(\A)) \\ \chi_R \simeq \pi_R}} \widehat{\phi}(\chi).
\end{equation}

Consider $\GL(2)$ instead of $\PGL(2)$ for simplicity, characters of $\PGL(2)$ corresponding to those of $\GL(2)$ trivial on the center.  Characters on $\GL(2)$ decompose through
\begin{equation}
G(F_\p) \longrightarrow F_\p^\times \longrightarrow S^1,
\end{equation}

\noindent where the first arrow is given by the determinant and the second by characters of $F_\p^\times$. In other words, a character $\chi_\p$ of $\GL(2, F_\p)$ is of the form $\chi_{0, \p} \circ \det$ where $\chi_{0, \p}$ is a character of $F_\p^\times$. 

At an archimedean place $v$, since the considered characters are trivial on the center, they are among the trivial one and the sign, hence have conductor 1 at those places. Archimedean characters are of the form $\mathrm{sgn}^{\varepsilon} |\det|^{it}$ for $\varepsilon = \pm 1$ and $t \in \R$. Similarly to the smoothing function introduced in Section \ref{subsec:arch-params}, it is not possible to select precisely continuous parameters, it is hence necessary to supply an approximation by a localizing function. This motivates the introduction of $f_v$ a compactly supported non-negative smooth function such that $\widehat{f}_v$ is 1 for $t=0$, and $|\widehat{f}_v| \leqslant 1$. In particular, it vanishes unless $t$ is small enough,  say $|t| \leqslant T$.

For the arithmetic part of the conductor, the only characters not killed by the action of $\widehat{\varepsilon}_{\p^r}$ are the unramified ones. Indeed, recall that
\begin{equation}
\det \left( K_{0}(\p^r)  \right) = \mathcal{O}_\p^\times,
\end{equation}

\noindent so that $\chi_{0, \p}$ needs to be trivial on $ \mathcal{O}_\p^\times$, that is to say be unramified. Introduce, for every finite split place $\p$, the characteristic function $f_\p$ of $\mathcal{O}_\p^\times$, whose Fourier transform selects unramified characters analogously to Lemma \ref{ftsapplied}. Introduce the global test function
\begin{equation}
f = \prod_{\p \notin R} f_{\p} \prod_{v \in R} \xi_{\chi_v} \prod_{\substack{v | \infty \\ v \notin R}} f_v.
\end{equation}

Since $\hat{f}_{\p}$ is 1 on unramified characters and the archimedean $\hat{f}_v$'s are less than one , the Poisson summation formula gives
\begin{equation}
\Xi (\pi_R, \phi) \ll \sum_{\chi \in \widehat{F^\times}} \widehat{f}(\chi) = \frac{1}{\vol(F^\times \backslash \A^\times)} \sum_{\gamma \in F^\times} f(\gamma).
\end{equation}

Since $F^\times$ is a discrete set, choosing $f_\infty$ with a small enough support leads to kill every $f(\gamma)$ for $\gamma$ nontrivial. Hence $\Xi (\pi_R, \phi) \leqslant \vol(F^\times \backslash \A^\times)^{-1} f(1)$. It remains to evaluate $f(1) = f_{\d^S}(1)f_R(1)f_{\infty} (1)$. For the finite split places, $f_{\p}(1) = 1$, and for the ramified places, $f_R(1) = \mu_R^{\mathrm{Pl}}(\pi_R)$. For the archimedean places, the Plancherel inversion formula gives
\begin{equation}
f_{\infty}(1) = \int_{\widehat{F}_\infty^R} \widehat{f}_{\infty}(\chi) \dd\chi \leqslant \int_{|t| \leqslant T} \dd \chi_t \ll_T 1.
\end{equation}

\noindent Finally, $\Xi(\pi_R, \PPP) \ll \mu_R^{\mathrm{Pl}}(\pi_R)$. Coming back to the sum \eqref{charac_contrib} defining $\nu_{\Xi, Q}(\widehat{\PPP}) $, it follows by using the rough bound $\lambda_2(\n) \ll N\n^\varepsilon$,
\begin{align*}
\nu_{\Xi, Q}(\widehat{\PPP}) & \ll \frac{1}{Q^2}\sum_{\substack{N\q \leqslant Q \\ \q \wedge R = 1}} \sum_{\substack{\pi_R \in \widehat{G}_R \\ c(\pi_R) \leqslant Q/N\q}} \sum_{\d^S \div \q^S} \lambda_2\left(\frac{\q^S}{\d^S}\right) \mu_R^\Pl(\pi_R) \\
& \ll \frac{1}{Q^2} \sum_{\substack{\pi_R \in \widehat{G}_R \\ c(\pi_R) \leqslant Q}} \mu_R^\Pl(\pi_R) \sum_{\substack{N\d \leqslant Q/c(\pi_R) \\ \q \wedge R = 1}} \sum_{\substack{N\m \leqslant Q/N\d c(\pi_R) \\ \q \wedge R = 1}} N\m^{\varepsilon} \\
& \ll Q^{-1+\varepsilon} \sum_{\substack{\pi_R \in \widehat{G}_R \\ c(\pi_R) \leqslant Q}} \frac{\mu_R^\Pl(\pi_R)}{c(\pi_R)^{1+\varepsilon}} \sum_{\substack{N\d \leqslant Q/c(\pi_R) \\ \q \wedge R = 1}} N\d^{-1-\varepsilon} \\
& \ll Q^{-1+\varepsilon}
\end{align*}

\noindent and this last line is provided by the convergence of the sum over ramified representation, stated in Lemma \ref{summeasures}, proving the result.  \qed

\subsection{Complementary spectrum}
\label{sec:complementary-spec}
\label{sec:COMP}

Lemma \ref{lem:selecting} provides an interpretation of $B(\q, \pi_R, \delta, Q, \phi)$ as the tempered part of the trace formula for a suitable test function, so it is necessary to consider the remaining complementary part of the spectrum. Lemma \ref{lem:arch-tf-spectral} states exponential control for spectral parameters lying outside the tempered subspace, so that the complementary contribution is bounded by $\partial_\rho B(\delta, \partial_\rho\Omega)$.  This section is dedicated to prove the following lemma, which is an adaptation of the work of Brumley and Mili\'{c}evi\'{c} in which the complementary spectrum is shown to contribute as an error term.

\begin{lemma}
\label{lem:BM-bounds-2}
The contribution of the complementary part and the smoothing error term satisfy, for a suitable choice of $\rho$ depending on $\q$ and a certain $\theta > 0$,
\begin{equation}
\sum_{\substack{N\q \leqslant Q \\ \q \wedge R, S = 1}}  \sum_{\substack{\sigma_R \in \widehat{G}_R \\ c(\sigma_R) \leqslant Q/N\q}} \sum_{\delta \in\mathcal{D}} \sum_{\d | \q} \lambda_2\left(\frac{\q^s}{\d^S} \right) \phi_2(\d^S)\mu_R^\Pl(\pi_R) B_{\mathrm{comp}} (\delta, Q(Q/N\q c(\pi_R))) \ll {Q^{2-\theta}} 
\end{equation}
\end{lemma}

\proof This is in essence Proposition 12.1 in \cite{brumley_counting_2018}. Their Lemma 12.2 yields
\begin{equation}
\sum_{\substack{N\q \leqslant Q \\ \q \wedge R, S = 1}} \sum_{\delta \in\mathcal{D}} N\q   B_{\mathrm{comp}} (\q, \delta, Q(Q/N\q)) \ll Q^{2-\theta} 
\end{equation}

Similarly to Lemma \ref{lem:BM-bounds}, the different congruence subgroups chosen for the centerless setting leads to input $N\q$ in the sum above, instead of $N\q^2$ in their case. The sum over ramified places is dealt with by appealing to Lemma \ref{summeasures} which ensures the convergence.
\begin{align*}
\sum_{\substack{N\q \leqslant Q \\ \q \wedge R, S = 1}} &  \sum_{\substack{\sigma_R \in \widehat{G}_R \\ c(\sigma_R) \leqslant Q/N\q}} \sum_{\delta \in\mathcal{D}} \sum_{\d | \q} \lambda_2\left(\frac{\q^S}{\d^S} \right) J_{\mathrm{comp}}\left(\Phi_{\d, \pi_R, \delta, Q_\delta(Q/N\d c(\pi_R))}\right)  \\
& \ll \sum_{\substack{N\q \leqslant Q \\ \q \wedge R, S = 1}}  \sum_{\substack{\sigma_R \in \widehat{G}_R \\ c(\sigma_R) \leqslant Q/N\q}} \sum_{\delta \in\mathcal{D}} \sum_{\d | \q} N\q \mu^\Pl_R(\pi_R) B_{\mathrm{comp}} (\delta, Q(Q/N\d c(\pi_R))) \\
& \ll Q^{2-\theta}\sum_{\substack{\sigma_R \in \widehat{G}_R \\ c(\sigma_R) \leqslant Q/N\q}} \frac{\mu^\Pl_R(\pi_R)}{c(\pi_R)^{2-\theta}}
\end{align*}

\noindent and this last sum converges by Lemma \ref{summeasures}, finishing the proof. \qed

\medskip

At last, it follows from this lemma along with Lemma \ref{lem:characrers} that the extra terms in the spectral selecting Lemma \ref{lem:selecting} are negligible, so that the counting measure part $B(\q, \pi_R, \delta, Q, \phi)$ is fairly well approximated by the tempered part $J_{\mathrm{temp}}(\Phi)$ of the spectrum. Lemma  \ref{lem:smoothing-part} states that this tempered part is approximated by the whole spectral part $J_{\mathrm{spec}}(\Phi)$ up to an error term. Considering the development of the identity contribution to the geometrical part given in Proposition \ref{prop:main-part-id-2}, it follows for every $\varepsilon >0$,
\begin{equation}
\nu_Q  = \vol\left(G\left(F\right) \backslash G\left(\A\right)\right) \nu_{1,Q} + \nu_{\ell, Q} + 
\left\{
\begin{array}{cl}
O\left( Q^{-1} \log(Q) \right) & \text{if $B$ totally definite and $F=\Q$;} \\
O\left( Q^{-\delta_F + \varepsilon}  \right) & \text{if $B$ totally definite and $F\neq\Q$;} \\
O\left( \log(Q)^{-1} \right) & \text{if $B$ not totally definite.} \\
\end{array}
\right.
\end{equation}

\section{Elliptic error terms}
\label{sec:elliptic}

The present section aims at bounding the different terms appearing in the elliptic contribution to the geometric side, in particular the orbital integrals. A considerable amount of work has been done in this direction, and we borrow recent results of \cite{binder_fields_2017}, \cite{finis_approximation_2013} and \cite{matz_sato-tate_2016} in order to reach our goal.

\subsection{Strategy}
\label{elliptic-strategy}

The contribution of the elliptic terms in the trace formula \eqref{fts} is
\begin{equation}
\label{fff}
J_{\ell}\left(\Phi\right) =  \sum_{\{\gamma\} \neq 1} \vol\left(G_\gamma(F) \backslash G_\gamma(\A)\right) \int_{G_\gamma(\A) \backslash G(\A)} \Phi\left(x^{-1}\gamma x\right) \dd x  .
\end{equation}

%
%


\noindent Recall that $\Phi$ denotes the test function $\Phi_{\d, \pi_R, \delta, Q, \rho; \PPP}$ introduced in Section \ref{subsec:fonction-test}, and that only a subset of the indices $(\d, \pi_R, \delta, Q, \rho, \PPP)$ may be used when the dependency on them has to be emphasized. As a matter of fact, the expression \eqref{fff} generally requires to bound
\begin{itemize}
\item the length of the summation, provided it is finite;
\item the global volumes $\mathrm{vol}\left(G_\gamma(F) \backslash G_\gamma(\A)\right)$;
\item the orbital integrals.
\end{itemize}

Since for every finite place $\p$, the test function $\Phi_\p$ is supported on either $K_\p$ in the case of a split place, or on $G_\p$ in the case of a ramified place, it is compactly supported on a compact independent of the fixed spectral parameters. However, the size of the sum is not uniformly bounded since the size of the support of the test function $f_\rho^{\delta, Q, \phi}$ at archimedean split places grows with the quality of the archimedean split approximation, encoded in $\rho$. Bounding the different quantities arising in the trace formula for archimedean split places has been addressed in many different works, and we heavily rely on \cite{brumley_counting_2018} and \cite{finis_approximation_2013}  to deal with these places.

\subsection{Contributing classes}

The conjugacy classes actually appearing in the sum are those for which $\mathcal{O}_\gamma(\Phi)$ does not vanish: let us call such conjugacy classes the contributing classes, according to \cite{matz_weyl_2017}. A contributing class therefore corresponds to a characteristic polynomial bounded to certain congruence and size conditions on its coefficients. This argument allows to count the number of conjugacy classes, and this is adapted to our setting in the following lemma \cite[Proposition 13.2]{brumley_counting_2018}.
\begin{lemma}
\label{lem:number}
For any split archimedean place $v$, there is a constant $c_v>0$ such that the number of conjugacy classes $\gamma_v$ of $G(F_v)$ for which $\mathcal{O}_{\gamma_v}(\Phi_{\rho, v})$ is nonzero is bounded by $\exp(c_v\rho)$.
\end{lemma}

Following the remark made above, the test function at any place that is not archimedean and split is uniformly compactly supported. Therefore, since there are only a finite number of split archimedean places, there is a constant $c > 0$ such that the number of contributing classes is bounded by $e^{c \rho}$.

\subsection{Global volumes}

The global volumes for $\GL(2)$ are bounded in \cite[Section 9]{matz_bounds_2015} and carefully adapted to our precise choice of archimedean test function in \cite[Proposition 12.4]{brumley_counting_2018}. We have the following.
\begin{lemma}
\label{lem:global-volumes}
There is a $c > 0$ such that for any $\rho > 0$ and for any $\gamma$ such that $\mathcal{O}_\gamma(\Phi_\rho)$ is nonzero, 
\begin{equation}
\vol\left(G_\gamma(F) \backslash G_\gamma(\A)\right) \ll e^{c \rho}.
\end{equation}
\end{lemma}

Finally, it remains to bound the orbital integrals defined by
\begin{equation}
\O_\gamma(\Phi) = \int_{G_\gamma(\A) \backslash G(\A)} \Phi(x\gamma x^{-1}) \dd x, \qquad \Phi \in \mathcal{H}(G_S).
\end{equation}

\noindent The needed bound is provided by the following proposition.

\begin{proposition}
\label{prop:iorb}
\label{prop:io}
For every $\gamma \in G(F)$, there is a $c > 0$ such that
\begin{equation}
\O_{\gamma}(\Phi) \ll_{\varepsilon} (N\d^R)^{-1+\varepsilon}  \mu^\Pl_R(\pi_R) \left\|f_\rho^{\delta, Q, \phi}\right\|_\infty e^{c \rho} .
\label{iorb}
\end{equation}
\end{proposition}

The next subsections are devoted to the proof of this proposition. The local components $\Phi_{\p}$ are almost always equal to $\mathbf{1}_{\overline{K}_{\p}}$, so that the corresponding local orbital integrals are almost always trivial by the normalizations of measures. Following \cite{knightly_traces_2006}, the local decomposition of orbital integrals for factorizable functions $\Phi = \otimes_v \Phi_v$ then holds, more precisely

\begin{equation}
\O_\gamma(\Phi) = \prod_v \O_{\gamma, v}(\Phi_v) \quad \text{where} \quad \O_{\gamma, v}(\Phi_v) = \int_{G_{\gamma, v} \backslash G_v} \Phi_v\left(x_v\gamma_v x_v^{-1}\right) \dd x_v.
\end{equation}

\noindent It thus suffices to dominate these local orbital integrals. The split, non-split and archimedean cases behave quite differently and require specific treatments. Also, it is critical to keep the dependence on $\gamma$ explicit and to formulate the bounds in terms of the Weyl discriminant of $\gamma$, heavily following the works of \cite{finis_approximation_2013}, \cite{matz_bounds_2015} and \cite{matz_sato-tate_2016}. The precise dependence on the size of the archimedean support, encoded in $\rho$, is also critical and the results aforementioned have been adapted to unveil this dependency by \cite{brumley_counting_2018}.

\subsection{Split orbital integrals}
\label{subsec:split-oi}

Almost every place is split, thus precise bounds are needed in order to control the global orbital integral. Fortunately, the test functions chosen at these places are explicit and allows to sharply control the associated orbital integrals.

\subsubsection{Non-archimedean split places}

\begin{lemma}
For every $\gamma^R \in G^R$, every ideal $\d^R$ of $\mathcal{O}^R$ and every $\varepsilon > 0$, 
\begin{equation}
\O_{\gamma^{R}}(\d^R) \ll \prod_{\p^r || \d^R} |D(\gamma)|_\p^{-1} N\left(\p^r\right)^{ \varepsilon}.
\end{equation}
\end{lemma}

\proof By the local factorization of orbital integral, it is sufficient to prove the lemma for a fixed place. Let $\p \notin R$ and $\gamma_\p \in G_\p$. In the case of a place $\p \notin S$, the local test function is of the form $\varepsilon_{\p^r}$, so that
\begin{equation}
\O_{\gamma}(\varepsilon_{\p^r}) = \vol(K)^{-1} \O_\gamma(\mathbf{1}_K), \qquad \text{ where } \quad \ K = \overline{K}_0(\p^r).
\end{equation}

\noindent Bounds for the split orbital integrals are provided by \cite[Proposition 8.2.1]{binder_fields_2017} in the specific case of $\GL(2)$, and yield the following estimate depending on $\gamma_\p$:
\begin{equation}
\label{binder-bound}
\O_{\gamma_\p}\left(\varepsilon_{\p^r}\right)  \ll  |D(\gamma)|_\p^{-1} N(\p^r)^{-1+\varepsilon} \vol(\overline{K}_0(\p^r))^{-1} \ll  |D(\gamma)|_\p^{-1} N\left(\p^r\right)^{\varepsilon}.
\end{equation}

\noindent  Otherwise, for $\p \in S$, the chosen test function is $\phi_\p$ and hence can be roughly bounded by
\begin{equation}
\mathcal{O}_{\gamma_\p} (\phi_\p) \ll |D(\gamma)|_\p^{-1},
\end{equation}

\noindent settling the desired estimates for finite split orbital integrals. \qed

\subsubsection{Archimedean split places}

\begin{lemma}
There is a constant $C > 0$ such that for every discrete spectral data $\delta \in \mathcal{D}$, every $\phi \in \mathcal{H}(G_\infty^R)$ and every $\rho, Q > 0$, we have
\begin{equation}
\O_{\gamma^R_\infty}\left( f_\rho^{\delta, Q, \phi} \right) \ll |D(\gamma)|^{-C} \left\| f_\rho^{\delta, Q} \right\|_\infty .
\end{equation}
\end{lemma}

\proof Using the fact that $\phi$ is bounded since compactly supported, the claim reduces to bound $\O_{\gamma^R_\infty}\left( f_\rho^{\delta, Q} \right)$. We can then appeal to \cite[Proposition 14.2]{brumley_counting_2018} which states the desired bound. \qed

\subsection{Non-split orbital integrals}
\label{subsec:nonsplit-io}

Ramified places are in finite number but the explicit behavior of local orbital integrals could a priori be unbounded. Underlining that the sum over conjugacy classes appearing as the geometric side of the trace formula is uniformly bounded, we can afford a dependence on $\gamma_v$ at ramified places. We have the following.
\begin{lemma}
\label{blou}
For every ramified place $v$,
\begin{equation}
\O_{\gamma_{v}}(\Phi_{v}) \ll |D(\gamma)|_v^{-1/2} \mu^\Pl_v(\pi_v).
\label{iorb}
\end{equation}
\end{lemma}

\proof Archimedean and non-archimedean ramified places behave differently. Before turning to the precise study of each case,  note that whatever $\Phi_v$ is $\xi_{\pi_v}$ or $\xi_{\pi_v} \widehat{\phi}_v(\pi_v)$, the orbital integral is dominated by the case of the matrix coefficient $\xi_{\pi_v}$, for $\widehat{\phi}_v$ is bounded. In the ramified case, orbital integrals are characters: for a representation $\pi_{v} \in \widehat{G}_v$, the main geometric lemma of \cite{arthur_characters_1988} implies that
\begin{equation}
\O_{\gamma_{v}}(\xi_{\pi_v}) \ll \Theta_{\pi_{v}}(\gamma_{v}) \mu^\Pl_v(\pi_v),
\label{xxxx}
\end{equation}

\noindent where $\Theta_{\pi_v}$ stands for the Harish-Chandra character associated to $\pi_v$.  It is in particular sufficient to bound characters on $B_\p^\times$ in order to get the desired bound for orbital integrals.

\subsubsection{Archimedean ramified places}

For matrix coefficients, we follow the work of \cite[Proposition 5.1]{kim_asymptotics_2016}. Since all the elements $\gamma_v$ are elliptic and regular, the Harish-Chandra character formula implies that 
\begin{equation}
\Theta_{\pi_{v}}(\gamma_{v}) \ll |D(\gamma)|_v^{-1/2}
\end{equation}

\subsubsection{Non-archimedean ramified places}

Concerning the non-archimedean ramified places $\p \in R$, the lead is given to \cite{shin_sato-tate_2016}, who build on the Sally-Shalika character formula in order to give explicit computations for the characters of each supercuspidal representations of $\SL(2)$. They prove that for every supercuspidal representation $\pi_\p$ of $\SL(2, F_{\p})$, and for all semisimple regular element $\gamma_\p$, 
\begin{equation}
\Theta_{\pi_\p}(\gamma_\p) \ll |D(\gamma)|_\p^{-1/2}.
\end{equation}

\noindent  Moreover, it suffices to achieve this goal for $\SL(2)$. Indeed, \cite{langlands_l-indistiguishability_1979} established that every irreducible admissible representation of $\GL(2)$ restricts to a direct sum of at most four irreducible admissible representations of $\SL(2)$. Since the Jacquet-Langlands correspondence maps irreducible representations of $G_\p$ to supercuspidal representations, and the image of the embedding of $G_\p$ in $\GL(2,F_\p)$ is made of semisimple regular elements, the bound above apply to $B_\p$, and therefore go $G_\p$.

\medskip

The bounds obtained in the two cases of ramified places hence settle the proof of Lemma \ref{blou}. \qed

\subsection{Sieving out Weyl's discriminants}

Let $\gamma \in G(F)$, $\d^R$ an ideal of $\mathcal{O}^R$, $\pi_R \in \widehat{G}_R$, $\delta \in \mathcal{D}$ and $\varepsilon, \rho, Q > 0$. Denote by $c_v > 0$ the values of the constants appearing in the lemmas above, which are uniformly bounded since they are almost always equal to $1/2$. Altogether, the previous local bounds compile into the global estimate
\begin{equation}
\mathcal{O}_\gamma(\Phi) \ll e^{c\rho} \left( N\d^S \right)^{-1+\varepsilon} \mu^\Pl_R(\sigma_R)\left\| f_\rho^{\delta, Q, \phi} \right\| \prod_v |D(\gamma)|_v^{-c_v}.
\end{equation}

It follows from \cite[Lemma 3.4]{matz_sato-tate_2016} that the archimedean norm of the Weyl discriminant is bounded as follows.
\begin{lemma}
There exists a constant $c > 0$ such that for every contributing class $\gamma \in G(F)$, $|D(\gamma)| \ll e^{c \rho}$.
\end{lemma}

For every finite place $\p$, let $\gamma_\p$ be the local component of a contributing class. Since the conjugacy class of $\gamma_\p$ meets $K_\p$, we have that $|D(\gamma)|_\p \leqslant 1$. We now argue as in the proof of the \cite{brumley_counting_2018} to get a uniform bound on $\gamma$. Since the contributing classes $\gamma$ are elements of $G(F)$, their Weyl discriminant is an element of $F$ and the product formula holds. Let $C > 0$ be the maximum of the $c_v > 0$, so that

\begin{align*}
\prod_v |D(\gamma)|_v^{-c_v}  & = \prod_v |D(\gamma)|_v^{C} \prod_v |D(\gamma)|_v^{-c_v} = \prod_v |D(\gamma)|_v^{C-c_v} \\
& \ll \prod_{v \, | \, \infty} |D(\gamma)|_v^{C-c_v} \ll e^{\kappa \rho}
\end{align*}

\noindent for a constant $\kappa > 0$ large enough. Altogether, we get the claimed bound in Proposition \ref{prop:io}. \qed

\subsection{Final estimates}
\label{sec:final}
\label{sec:final-estimates}

We first recall the lemma \cite[Lemma 7.3]{brumley_counting_2018} that contains the behavior of the sum over the discrete spectral data of the bounds found in the previous paragraphs. 
\begin{lemma}
\label{lem:archIO}
For every $X>0$, there is $c > 0$ so that
\begin{equation}
\label{bound-arch-OI}
\sum_{\delta \in \mathcal{D}} \left\|f_\rho^{\delta, Q_{\delta}(X), \phi }\right\|_\infty \ll e^{c \rho} X^{2-1/[F:\Q]}.
\end{equation}
\end{lemma}

\noindent All the tools are now in place to establish the final estimates on the elliptic contribution $\nu_{Q, \ell}(\widehat{\phi}\,)$ and reach the term of the proof of Theorem \ref{thm:equid-prec}.
\begin{proposition}
\label{elliptic}
Let $d = [F:\Q]$. For a finite set of places $S$, $\phi \in \mathcal{H}(G_S)$ and any $\varepsilon > 0$, the elliptic contribution is dominated  by
\begin{equation}
\nu_{\ell, Q}(\widehat{\phi}\,) \ll Q^{-1/d + \varepsilon}.
\end{equation}
\end{proposition}

\proof The bounds stated in Proposition \ref{prop:io} and the definition of the elliptic contribution lead to
\begin{align*}
J_{\ell}\left(\Phi_{\d,\pi_R, \delta, Q;  \PPP}\right) & = \sum_{\{\gamma\} \neq 1} \vol(\Gamma_\gamma \backslash G_\gamma) \O_\gamma\left(\Phi_{\d,\pi_R, \delta, Q;  \PPP}\right) \\
& \ll  (N\d^S)^{\varepsilon} \mu^\Pl_R(\pi_R) e^{c \rho}  \left\|f_\rho^{\delta, \Omega_{\delta}(Q/N\q^S c(\pi_R)), \phi }\right\|_\infty 
\end{align*}

\noindent Summing over the spectral data and using Lemma \ref{lem:archIO} give
\begin{align*}
\nu_{\ell, Q}(\widehat{\phi}\,) & \ll \frac{1}{Q^2} \sum_{\substack{N\q \leqslant Q \\ \q \wedge R = 1}} \sum_{\d^S\div \q^S} \lambda_2\left(\frac{\q^S}{\d^S}\right) (N\d^S)^{\varepsilon}\sum_{\substack{\pi_R \in \widehat{G}_R \\ c(\pi_R) \leqslant Q/N\q^S}}  \widehat{\PPP}(\pi_R) \mu^\Pl_R(\pi_R) \sum_{\substack{\underline{\delta} \in \mathcal{D} \\ \underline{\delta} = (M, \delta)}} e^{c \rho} \left\|f_\rho^{\delta, \Omega_{\delta}(Q/N\q^S c(\pi_R)), \phi }\right\|_\infty \\
& \ll  Q^{-1/d}  \sum_{\substack{N\q \leqslant Q \\ \q \wedge R = 1}} \sum_{\d^S\div \q^S} \lambda_2\left(\frac{\q^S}{\d^S}\right) e^{c \rho} (N\d^S)^{\varepsilon+1/d-2}\sum_{\substack{\pi_R \in \widehat{G}_R \\ c(\pi_R) \leqslant Q/N\q^S}}  \frac{\widehat{\PPP}(\pi_R) \mu^\Pl_R(\pi_R)}{c(\pi_R)^{2-1/d}}
\end{align*}

\noindent where the bound \eqref{bound-arch-OI} and the elementary bound $\lambda_2(\mathfrak{n}) \ll_\varepsilon N\mathfrak{n}^\varepsilon$ have been used. Thus, since Lemma \ref{summeasures} and \ref{summeasures2} ensure the convergences of the inner sum, it follows that for $\rho = \alpha \log N\d^S$ where $\alpha > 0$ is small enough,
\begin{align*}
& \nu_{\ell, Q}(\widehat{\phi}\,)   \ll_\varepsilon  Q^{-1/d+\varepsilon}, 
\end{align*}

\noindent and this choice of $\rho$ is admissible by the work of Brumley and Mili\'{c}evi\'{c}. This achieves the proof that the main term contributing in \eqref{splitting} is the one coming from the identity as stated in Proposition \ref{prop:id}, hence also Theorems \ref{thm-count}, \ref{thm-equid} and \ref{thm:equid-prec}.  \qed

\section{Sato-Tate corollary}
\label{sec:ST}

Theorem \ref{thm-equid} proves the existence of a measure $\nu$ with respect to which the universal family equidistributes. Consider the projection $\nu_\p$ of $\nu$ on the local components $\widehat{G}_\p$. Since the $\nu_\p$ are supported on different spaces, it is necessary to make sense of the Sato-Tate\index{Sato-Tate} problem that concerns convergence of the measures $\nu_\p$. 

The literature often treat the case of measures supported on the unramified tempered spectrum, as the instances handled by \cite{sarnak_statistical_1987} or \cite{serre_repartition_1997}. In those cases, the Satake isomorphism provides a common parametrization: if $T$ is the standard torus of $\SL(2, \C)$, the dual group of $\PGL(2)$, and $W$ is the Weyl group\index{Weyl group, $W$} of $T$, then the isomorphism classes of unramified tempered representations are parametrized by $T_c/W$ where $T_c = T \cap SU(2, \C)$ is the compact part of $T$. This last quotient corresponds to the half-circle, giving a common ground for all the $\widehat{G}_\p$,independent of $\p$. Even if the universal family considered does include ramified representations and the $\nu_\p$ are supported on the whole tempered unitary dual, the contribution of the ramified part of the spectrum vanish when $\p$ goes to infinity, so that asymptotically the spaces can be identified and $T_c/W$ is a posteriori still a relevant common ground to state the Sato-Tate result.

For $\GL(2,F_\p)$, the Plancherel measures have been computed by \cite{serre_repartition_1997} and are given by
\begin{equation}
\label{plancherel-measure-local}
\dd\mu^\Pl_\p(x)  = \frac{N\p+1}{\pi} \frac{(1-x^2/4)^{1/2}}{(N\p^{1/2}+N\p^{-1/2})^2 - x^2} \dd x.
\end{equation}

\noindent In particular they converge, as $N\p$ goes to infinity, to the Sato-Tate measure on the half-circle
\begin{equation}
\dd\mu^{\mathrm{ST}}(x) = \frac{1}{\pi} \sqrt{1-\frac{x^2}{4}}\dd x,
\end{equation}

\noindent in the sense that for any $\widehat{\phi} \in C(T_c/W, \C)$, when $N\p$ goes to infinity,
\begin{equation}
\int_{T_c/W}\widehat{\phi}(\pi_\p)\dd\mu^\Pl_\p (\pi_\p) \longrightarrow \int_{T_c/W}\widehat{\phi}(x)\dd\mu^{\mathrm{ST}}(x).
\end{equation}


For $\widehat{\phi} \in C(T/W, \C)$, let decompose the measure separating whether the representations are unramified, \textit{i.e.} of conductor 1, or not. The measure $\nu_\p(\widehat{\phi}\,)$ hence splits as
\begin{equation}
\begin{split}
\int_{\widehat{G}_\p}\widehat{\phi}(\pi_\p)\dd\nu_\p(\pi_\p) & = \int_{\widehat{G}_\p}\frac{\widehat{\phi}(\pi_\p)}{c(\pi_\p)^2}\dd\mu^\Pl_\p (\pi_\p) \\
& = \int_{\widehat{G}_\p^\mathrm{sph}}\widehat{\phi}(\pi_\p)\dd\mu^\Pl_\p (\pi_\p) + \int_{\widehat{G}_\p^\mathrm{ram}} \frac{\widehat{\phi}(\pi_\p)}{c(\pi_\p)^2}\dd\mu^\Pl_\p (\pi_\p),
\end{split}
\end{equation}

\noindent where $\widehat{G}_\p^\mathrm{sph}$ stands for the unramified, also called spherical, part of the spectrum and $\widehat{G}_\p^\mathrm{ram}$ for its ramified part. For $\p$ sufficiently large, $G_\p$ is isomorphic to $\PGL(2, F_\p)$, so the local Plancherel measures \eqref{plancherel-measure-local} provide the value of the first integral of the rightmost hand side as $\p$ grows, in particular they converge to the Sato-Tate measure. For the second one, dominating roughly by leaving the dependence in $\phi$ which is fixed gives
\begin{equation}
\int_{\widehat{G}_\p^\mathrm{ram}} \frac{\widehat{\phi}(\pi_\p)}{c(\pi_\p)^2}\dd\mu^\Pl_\p (\pi_\p) \ll \int_{\widehat{G}_\p^\mathrm{ram}}\frac{\dd\mu^\Pl_\p (\pi_\p)}{c(\pi_\p)^2} = \int_{\widehat{G}_\p}\frac{\dd\mu^\Pl_\p (\pi_\p)}{c(\pi_\p)^2} - \int_{\widehat{G}_\p^\mathrm{sph}}\dd\mu^\Pl_\p (\pi_\p).
\end{equation}

By the normalization of the Plancherel measure, the second integral on the right hand side is 1. Moreover, as shown in Section \ref{subsec:constant}, the first integral of the right hand side is equal to
\begin{equation}
\int_{\widehat{G}_\p}\frac{\dd\mu^\Pl_\p (\pi_\p)}{c(\pi_\p)^2} = \frac{\zeta_\p(1)}{\zeta_\p(2)\zeta_\p(4)}.
\end{equation}

Since this last quantity is $1+O(N\p^{-1})$ by unfolding the definition of the Dirichlet series, it follows that the ramified part is negligible, achieving the proof of Corollary \ref{coro:ST}. \qed


\begin{thebibliography}{00}

\bibitem[{Adler et al.},{2011}]{adler_supercuspidal_2011}
	{Adler, J. D. and DeBacker, S. and Sally, P. Jr. and Spice, L.}
	\textit{Supercuspidal characters of {SL}(2) over a p-adic field}.
	In: {Harmonic {Analysis} on {Reductive}, p-adic {Groups}},
	pp. {19--69}.
	{Amer. Math. Soc.}
	(2011)
	
	\bibitem[{Arthur},{1988}]{arthur_characters_1988}
	{Arthur, J.}
	{The characters of supercuspidal representations as weighted orbital integrals}.
	{Proc. Indian Acad. Sci. Math. Sci}.
	{97},
	{3--19}
	({1988})
	
\bibitem[{Arthur},{2005}]{arthur}
	{Arthur, J.}
	\textit{An introduction to the trace formula}.
	In:Harmonic analysis, the trace formula, and {Shimura} varieties.
	{Clay. {Math}. {Proc}.} 4,
	{1--263}
	(2005)
	
\bibitem[{Batyrev et al.},{1990}]{batyrev_sur_1990}
	{Batyrev, V. V. and Manin, Y. I.}
	\textit{Sur le nombre des points rationnels de hauteur bornée des variétés algébriques}.
	 {Math. Ann.}.
	 286,
	 {27- 43}
	(1990)
	
	\bibitem[{Bernstein et al.},{1986}]{bernstein_trace_1986}
	{Bernstein, J. and Deligne, P. and Kazhdan, D.}
	\textit{Trace paley-wiener theorem for reductivep-adic groups}.
	 {J. Analyse Math}.
	 47,
	{180--192}
	(1986)
	
	\bibitem[{Binder},{2017}]{binder_fields_2017}
	{Binder, J.}
	\textit{Fields of {Rationality} of {Cusp} {Forms}}.
	 {Israel J. Math}.
	 222,
	{973--1028}
	(2017)
	
\bibitem[{Brooks et al.},{2018}]{brooks_counting_2016}
	{Brooks, E. H. and Petrow, I.}
	\textit{Counting {automorphic} {forms} on {norm} {one} {tori}}
	{Acta Arith.}.
	183, 
	{117--143}
	(2018)
	
	
\bibitem[{Brumley},{2006}]{brumley_effective_2006}
	{Brumley, F.}
	\textit{Effective {Multiplicity} {One} for {GL}({N}) and narrow zero-free regions for {Rankin}-{Selberg} {L}-functions}.
	 {Amer. J. Math}.
	 {128},
	{1455--1474}
	(2006)
	
			\bibitem[{Brumley et al.},{2018}]{brumley_counting_2018}
    {Brumley, F. and Mili\'{c}evi\'{c}, D.}
	\textit{Counting cusp forms with analytic conductor}
	(2018).
	Preprint {\href{https://arxiv.org/pdf/1805.00633.pdf}{arXiv:1805.00633}}
	
	\bibitem[{Casselman},{1973}]{casselman_results_1973}
	{Casselman, W.}
	\textit{On {Some} {Results} by {Atkin} and {Lehner}}.
	{Math. Ann}.
	{201},
	 {301--314}
	(1973)
	
\bibitem[{Chamnert-Loir et al.},{2010}]{chambert-loir_igusa_2010}
	{Chambert-Loir, A. and Tschinkel, Y.}
	\textit{Igusa integrals and volume asymptotics in analytic and adelic geometry}.
	{Confluentes Math}.
	{2},
	{351--429}
	(2010)
	
		\bibitem[{Clozel et al.},{1990}]{clozel_theoreme_1990}
	{Clozel, L. and Delorme, P.}
	\textit{Le théorème de {Paley}-{Wiener} invariant pour les groupes réductifs. {II}}.
	{Ann. Sci. ÉNS}.
	{23},
	{193--228}
	(1990)
	
	\bibitem[{Conrey et al.},{2005}]{conrey_integral_2005}
	{Conrey, J. B. and Farmer, D. W. and Keating, J. P. and Rubinstein, M. O. and Snaith, N. C.}
	\textit{Integral moments of {L}-functions}.
	{Proc. London Math. Soc.}.
	{91},
	{33--104}
	(2005)
	
\bibitem[{Diamond et al.},{2005}]{diamond_first_2005}
	{Diamond, F. and Shurman, J.}
	\textit{A {First} {Course} in {Modular} {Forms}}.
	{Graduate {Texts} in {Math}} 228.
	Springer
	(2005)
	
		\bibitem[{Duistermaat et al.},{1979}]{DKV}
	{Duistermaat, J. J. and Kolk, J. A. C. and Varadarajan, V. S.}
	\textit{Spectra of compact locally symmetric manifolds of negative curvature}.
	{Inventiones mathematicae}.
	{52},
	{27--93}
	(1979)
	
	\bibitem[{Finis et al.},{2013}]{finis_approximation_2013}
	Finis, T. and Lapid, E.
	\textit{An approximation principle for congruence subgroups}.
	(2013)
	Preprint {\href{https://arxiv.org/pdf/1808.09991.pdf}{arXiv:1808.09991}}
	
\bibitem[{Godement et al.},{1972}]{godement_zeta_1972}
	{Godement, R. and Jacquet, H.}
	\textit{Zeta {Functions} of {Simple} {Algebras}}.
	Lecture Notes in Math 260.
	Springer 
	(1972)
	
		\bibitem[{Hare},{1998}]{hare_size_1998}
	{Hare, K. E.}
	\textit{The size of characters of compact {Lie} groups}.
	{Studia Math.}.
	{129},
	1--18
	(1998)
	
	\bibitem[{Hull},{1936}]{hull}
	{Hull, R.}
	\textit{The Maximal Order of Generalized Quaternion Division Algebras}.
	{Bull. Amer. Math. Soc.}.
	22, 
	1--11
	(1936)
	
	\bibitem[{Iwaniec et al.},{2000}]{iwaniec_perspectives_2000}
	{Iwaniec, H. and Sarnak, P.}
	\textit{Perspectives on the {Analytic} {Theory} of {L}-functions}.
	 {Geom. Funct. Anal}.
	{705--741},
	(2000)
	
		\bibitem[{Kim et al.},{2016}]{kim_asymptotics_2016}
	{Kim, J.-L., Shin, S. W. and Templier, N.}
	\textit{Asymptotics and Local Constancy of Characters of p-Adic Groups}.
	In: {Simons Symposia},
	{259–295}.
	Springer
	(2016)
	
		\bibitem[{Knapp},{1986}]{knapp_representation_1986}
	{Knapp, A. W.}
	\textit{Representation {Theory} of {Semisimple} {Groups}: {An} {Overview} based on {Examples}}.
	Princeton {Landmarks} in {Math} 36
	Princeton University Press
	(1986)
	
	\bibitem[{Knapp},{1994}]{knapp_local_1994}
{Knapp, A. W.}
	\textit{Local {Langlands} {Correspondence} {The} {Archimedean} {Case}}.
	In: {Proc. {Symp}. {Pure} {Math}} 55,
	{393--410}.
	Amer. Math. Soc.
	(1994)
	
	\bibitem[{Knapp et al.},{2000}]{knapp_representation_2000}
	{Knapp, A. W. and Trappa, P. E.}
	\textit{Representation of {Semisimple} {Lie} {Groups}}.
	In: {Representation {Theory} of {Lie} {Groups}},
	pp. {5--88}.
	{Amer. Math. Soc}
	(2000)

\bibitem[{Knightly and al.},{2006}]{knightly_traces_2006}
	{Knightly, A. and Li, C.}
	\textit{Traces of {Hecke} operators}.
	{Math. {Surveys} and {Monographs}} {133}.
	{Amer. Math. Soc}
	(2006)
	
	\bibitem[{Landau},{1918}]{landau_einfuhrung_1918}
	{Landau, E.}
	\textit{Einführung in die elementare und analytische {Theorie} der algebraischen {Zahlen} und der {Ideale}}.
	{Leipzig, B. G. Teubner}
	(1918)
	
	\bibitem[{Lang},{1985}]{lang_sl2_1985}
	{Lang, S.}
	\textit{SL(2,R)},
	{Graduate {Texts} in {Mathematics}} 105,
	Springer-Verlag,
	(1985)

\bibitem[{Langlands et al.},{1979}]{langlands_l-indistiguishability_1979}
	{Langlands, R. P. and Labesse, J.-P.}
	\textit{L-indistinguishability for {SL}(2)}.
	{Can. J. Math}.
	31,
	 {726--785}
	(1979)
	
	\bibitem[{Lansky et al.},{2004}]{lansky}
	{Lansky, J. and Raghuram, A.}
	\textit{On the correspondence of representations between {GL}(n) and division algebras}.
	 {Proc. Amer. Math. Soc.}.
	 131,
	 1641--1648
	(2004)
	
	\bibitem[{Matz},{2015}]{matz_bounds_2015}
	{Matz, J.}
	\textit{Bounds for global coefficients in the fine geometric expansion of Arthur's trace formula for GL(n)}.
	 {Israel J. Math.}.
	 205,
	 337--396
	(2015)
	
	\bibitem[{Matz},{2017}]{matz_weyl_2017}
	{Matz, J.}
	\textit{Weyl's law for Hecke operators on GL(n) over imaginary quadratic number fields}.
	 {Amer. J. Math.}.
	 139,
	 57--145
	(2017)
	
		\bibitem[{Matz et al.},{2016}]{matz_sato-tate_2016}
	Matz, J. and Templier, N.
	\textit{Sato-{Tate} equidistribution for families of Hecke-Maass forms on SL(n, R)/SO(n)}.
	(2016)
	Preprint {{arXiv:1808.09991}}
	
	\bibitem[{Northcott},{1950}]{northcott_periodic_1950}
	{Northcott, D. G.}
	\textit{Periodic {Points} on an {Algebraic} {Variety}}.
	{Ann. of Math}.
	{51},
	 {167--177}
	(1950)
	
\bibitem[{Petrow},{2018}]{petrow_weyl_2018}
	{Petrow, I.}
	\textit{The {Weyl} law on algebraic tori}
	(2018).
	Preprint {\href{https://arxiv.org/pdf/1808.09991.pdf}{arXiv:1808.09991}}
	
		
		\bibitem[{Peyre},{1995}]{peyre_hauteurs_1995}
	{Peyre, E.}
	\textit{Hauteurs et mesures de {Tamagawa} sur les variétés de {Fano}}.
	{Duke J. Math.}.
	{79},
	{101--218}
	 (1995)
	 
	
\bibitem[{Prasad et al.},{2000}]{prasad_representation_2000}
	{Prasad, D. and Raghuram, A.}
	\textit{Representation {Theory} of {GL}(n) over {Non}-{Archimedean} {Local} {Fields}}.
	{Lectures given at the School on Automorphic Forms on GL(n)}.
	{1-48} 
	(2000)
	
\bibitem[{Sally et al.},{1968}]{sally_characters_1968}
	{Sally, P. J. and Shalika, J.}
	\textit{Characters of the {Discrete} {Series} of {Representations} of {SL}(2) over a {Local} {Field}}.
	 {Proc. National Academy of Sciences of the United States of America}.
	 {61},
	{1231--1237}
	(1968)
	

\bibitem[{Sarnak},{1987}]{sarnak_statistical_1987}
	{Sarnak, P.}
	\textit{Statistical {Properties} of {Eigenvalues} of the {Hecke} {Operators}}.
	In: {Analytic {Number} {Theory} and {Diophantine} {Problems}},
   {321--331}.
   {Birkhäuser}
	(1987)

\bibitem[{Sarnak et al.},{2016}]{sarnak_families_2016}
	{Sarnak, P. and Shin, S. W. and Templier, N.}
	\textit{Families of {L}-functions and their symmetries}.
	In: {Families of automorphic forms and the trace formula},
   {531--578}.
   {Springer}
	(2016)

\bibitem[{Shin et al.},{2016}]{ST}
Shin, S. W. and Templier, N.
\textit{Sato--Tate theorem for families and low-lying zeros of automorphic                                                           L-functions.}
Inventiones mathematicae.
203,
 1--177
(2016)

\bibitem[{Sauvageot},{1997}]{sauvageot_principe_1997}
	{Sauvageot, F.}
	\textit{Principe de densit\'{e} pour les groupes réductifs}.
	{Compositio Math}.
	{108},
	{151--184}
	(1997)

\bibitem[{Schanuel},{1964}]{schanuel_heights_1964}
	{Schanuel, S.}
	\textit{On heights in number fields}.
	{Bull. Amer. Math. Soc}.
	{70},
	 {262--263}
	 (1964)

\bibitem[{Serre},{1997}]{serre_repartition_1997}
	{Serre, J.-P.}
	\textit{Répartition asymptotique des valeurs propres de l'opérateur de {Hecke} {T}(p)}.
	{J. Amer. Math. Soc}.
	{10},
	{75--102}
	(1997)
	
\bibitem[{Shin et al.},{2016}]{shin_sato-tate_2016}
	Shin, S. W. and Templier, N.
	\textit{Sato-{Tate} theorem for families and low-lying zeros of automorphic {L}-functions}.
	{Invent. Math}.
	{1}, 
	{1--177}
	(2016)
\end{thebibliography}
\end{document}